\documentclass[a4paper,12pt]{article}
\usepackage{mathrsfs}

\usepackage{latexsym}

\usepackage[mathcal]{euscript}

\usepackage[english]{babel}

\usepackage[ansinew]{inputenc}

\usepackage{amsfonts}
\usepackage{amsmath}
\usepackage{amssymb}

\usepackage[all]{xypic}

\long\def\symbolfootnote[#1]#2{\begingroup%
\def\thefootnote{\fnsymbol{footnote}}\footnote[#1]{#2}\endgroup} 

\newcommand{\bysame}{\leavevmode\vrule height 2pt depth -1.6pt width 23pt \ }
\newcommand{\mfdm}{\mathcal{M}}
\newcommand{\kd}{d_{\mathcal{K}}}
\newcommand{\dw}{d_{\mathcal{W}}}
\newcommand{\tncwa}{\mathcal{W}_{\g}^{(\chi)}}
\newcommand{\ncf}{C^{\infty}(\mathcal{M}_{\theta})}

\newcommand{\hmod}{_{\mathcal{H}}\mathscr{M}}
\newcommand{\h}{\mathfrak{h}}
\newcommand{\p}{\mathfrak{p}}
\newcommand{\pt}{\tilde{\p}}
\newcommand{\kk}{\mathfrak{k}}
\newcommand{\kkt}{\tilde{\kk}}
\newcommand{\g}{\mathfrak{g}}
\newcommand{\ttt}{\mathfrak{t}}
\newcommand{\sg}{\bar{\g}}
\newcommand{\gt}{\tilde{\mathfrak{g}}}
\newcommand{\gtb}{\sg^B}
\newcommand{\nwd}{d_{\mathcal{W}}}

\newcommand{\env}[1]{\mathfrak{U}(#1)}
\newcommand{\envt}[1]{\mathfrak{U}^{\chi}(#1)}
\newcommand{\ncwa}{\mathcal{W}_{\g}}
\newcommand{\proof}[1]{\noindent\underline{\textsf{Proof}}: #1 \hspace{0.5cm} $_\blacksquare$\\}

\newcommand{\car}[1]{C_G(#1)}
\newcommand{\nccarl}[1]{(\env{\g}\otimes #1)^G}
\newcommand{\nccar}[1]{\mathcal{C}_G(#1)}
\newcommand{\subnccar}[2]{\mathcal{C}_{#1}(#2)}
\newcommand{\tnccar}[1]{(\{K_a\}\otimes #1_{\chi})^G}
\newcommand{\tnccardf}[1]{\mathcal{C}^{\chi}_G(#1_{\chi})}
\newcommand{\subtnccardf}[2]{\mathcal{C}^{\chi}_{#1}(#2_{\chi})}

\newcommand{\cwm}{\left((\ncwa\otimes\mathcal{A})_{bas}\right)_{[[\theta]]}}
\newcommand{\twm}{(\tncwa\widehat{\otimes}\mathcal{A}_{\chi})_{bas}}

\newtheorem{df}{Definition}[section]
\newtheorem{ex}[df]{Example}
\newtheorem{thm}[df]{Theorem}
\newtheorem{prop}[df]{Proposition}
\newtheorem{lemma}[df]{Lemma}

\begin{document}

\begin{flushright}
SISSA Preprint $42/2007/$MP
\end{flushright}


\begin{center}
    \textbf{\Large{Twisted noncommutative  equivariant cohomology:   Weil and Cartan models}}
\vspace{0.4cm}

\textmd{Lucio Cirio\symbolfootnote[2]{Current address: Max Planck Institute for Mathematics - Vivatsgasse 7,
53111 Bonn, Germany. Email: cirio@mpim-bonn.mpg.de} \\ \vspace{0.25cm} International School for Advances Studies \\ Via Beirut $2-4$, $34014$ Trieste, Italy}

\vspace{0.4cm}
\textbf{Abstract}

\end{center}

\noindent
We propose Weil and Cartan models for the equivariant cohomology of noncommutative spaces which carry a covariant action of Drinfel'd twisted symmetries. The construction is suggested by the noncommutative Weil algebra of Alekseev and Meinrenken \cite{am1}; we show that one can implement a Drinfel'd twist of
their models in order to take into account the noncommutativity of the spaces we are acting on.

\tableofcontents

\clearpage

\section*{Introduction}

\addcontentsline{toc}{section}{Introduction}

The main goal of this paper is to introduce algebraic models for the equivariant cohomology of noncommutative (we use the shorter 'nc' throughout the paper) spaces acted covariantly by symmetries deformed by Drinfel'd twists. The covariance between the symmetries acting and the spaces acted is expressed by working in the category of Hopf-module algebras; we are eventually interested in defining equivariant cohomology of deformed Hopf-module algebras of a certain kind. 

More in detail, we start by considering actions of compact Lie groups $G$ on smooth manifolds $\mfdm$. The equivariant cohomology ring $H_G(\mfdm)$ is a natural tool for the study such actions. It replaces of the ordinary cohomology ring of the space of orbits $H(\mfdm/G)$ when the latter is not defined; the topological Borel model computes $H_G(\mfdm)$ as the ordinary cohomology of $EG\times_G\mfdm$, where $EG$ is the total space of the universal $G$-bundle. It is often convenient to switch to an algebraic description of the Borel model, replacing the infinite dimensional space $EG$ by a finitely generated algebra representing its differential forms, the Weil algebra $W_{\g}=Sym(\g^{\ast})\otimes\wedge(\g^{\ast})$. In this way we obtain the Weil model for equivariant cohomology, defined as the cohomology of the basic subcomplex of $W_{\g}\otimes\Omega(\mfdm)$. Another equivalent algebraic definition of $H_G(\mfdm)$, closer to the definition of de Rham cohomology of $\mfdm$, is formulated by introducing equivariant differential forms and then taking cohomology with respect to an equivariant differential operator $d_G$; this is known as the Cartan model. Ax excellent review on these classical models is \cite{gs}.

Both Weil and Cartan models make a crucial use of the operators $(i,L,d)$ (respectively interior, Lie and exerior derivative) on $\Omega(\mfdm)$. This triple provides a purely algebraic description of the action, as first observed by Cartan \cite{car}. It is convenient to introduce a super Lie algebra $\gt$, whose Lie brakets are given by the commutation relations of $(i,L,d)$; then $\Omega(\mfdm)$ carries a representation of $\gt$ by graded derivations, or equivalently a $\env{\gt}$-module structure. In the spirit of nc geometry, thus referring to algebras rather than spaces, we may say that classical equivariant cohomology is defined on the cateogry of $\env{\gt}$-module algebras. The Hopf module structure gives a link between deformations of the space (the algebra $\Omega(\mfdm)$) and deformation of the symmetry (the Hopf algebra $\env{\gt}$). We focus on deformations described by Drinfel'd twists of $\env{\gt}$, and we address the problem to define equivariant cohomology for twisted $\env{\gt}$-module algebras; this description applies to a quite general class of nc spaces, notably toric isospectral deformations.

The definition of a Weil model for deformed $\env{\gt}$-module algebras goes through a deformation of the classical Weil algebra $W_{\g}$. We take inspiration by the work of Alekseev and Meinrenken and their nc Weil algebra $\ncwa$ \cite{am1}\cite{am2}; we remark however that their models apply to nc $\env{\gt}$-module algebras, so with a classical action of $(i,L,d)$, while we are interested in algebras where the noncommutativity is strictly related to the deformation of the $\env{\gt}$-module structure. Our idea is that depending on the category of algebras we are interested in, the universal locally free object may be interpreted as the relevant Weil algebra for the definition of an equivariant cohomology. The nc Weil algebra $\ncwa$ comes by considering the category of nc $\env{\gt}$-module algebras; we are interested in the correspondind Drinfel'd twisted category, and we show indeed that it is possible to implement a twist of the construction of \cite{am1} to adapt the models to the class of nc spaces we study. We present also some examples of this twisted nc equivariant cohomology, and we discuss the property of maximal torus reduction $H_G(\mfdm)\cong H_T(\mfdm)^W$. We make some comments on the quite classical behavior of the defined cohomology, and we conclude by sketching how it could be possible to define equivariant cohomology of $\env{\gt}$-module algebras subjected to more general classes of deformations. 

This paper is structured in two sections. Section $1$ focuses on the relation between deformed symmetries and nc spaces; we start by motivating our interest on the category of Hopf-module algebras and we discuss some properties which will be relevant for the rest of the paper, notably its braided structure. We then introduce Drinfel'd twists on Hopf algebras and we characterize their effect on the category of modules, and as an example of nc spaces 'generated' by Drinfel'd twists we review toric isospectral deformations in this formalism. Section $2$ is devoted to equivariant cohomology; we start by recalling the classical definitions and models, then we move to the nc equivariant cohomology of Alekseev and Meinrenken. We finally introduce our twisted models, providing examples and discussing the reduction to the maximal torus. The section ends with a summary of a five-steps strategy towards the definition of equivariant cohomology for further classes of nc spaces.\\

\noindent
\textbf{Acknowledgments}\\

\noindent This work is part of my PhD thesis; I would like first if all to thank my advisors U. Bruzzo and G. Landi for many helpful suggestions and remarks. I wish to thank also S. Brain, F. D'Andrea and E. Petracci for interesting discussions during the various stages of this paper.

\section[Covariant actions on noncommutative spaces]{Covariant actions on noncommutative spaces}

In this section we describe symmetries of nc spaces; our approach focuses on the link between a symmetry and the space acted on by using the language of Hopf module algebras. As we will show one can deform a symmetry (as Hopf algebra) and then induce a compatible deformation in every space acted on. On the other hand one can rather start with a nc space and deform its classical symmetries in a suitable way in order to have a covariant action (see definition below) in the nc setting. Broadly speaking, while in the first case we 'adapt' classical spaces to deformed symmetries, thus regarding the symmetry as the 'source' of the nc deformation, in the second situation we force symmetries to fit a 'pre-existing' noncommutativity of spaces. 

In the first subsection we explain why the category of Hopf module algebras provides the natural setting where to study actions of classical and deformed symmetries; we also introduce some notions of braided categories which turn out to be very useful to describe nc spaces as 'braided-commutative' spaces. In the second subsection we focus on a particular class of Hopf algebra deformations, namely the ones coming from a Drinfel'd twist; we recall definitions and general properties of such twisted Hopf algebras. In the third subsection we discuss an interesting and well known example of nc spaces obtained from a twisted symmetry: toric isospectral deformations; we also show how to deform further symmetries acting on twisted nc spaces in order to preserve the covariance of the action. This will be used in Section 2 to define algebraic models for equivariant cohomology of such deformed actions.

\subsection{Hopf-module algebras}

We begin by describing the action of a compact Lie groups $G$ on a smooth compact Hausdorff manifold $\mfdm$ into a purely algebraic formalism. This language was introduced by H. Cartan \cite{car}, and it belongs by now to a classical background of differential geometry; for a modern and detailed treatment an excellent reference is \cite{gs}.

Let $A=\Omega^{\bullet}(\mfdm)$ be the graded-commutative algebra of differential forms
on $\mfdm$, and $\g$ the Lie algebra of $G$ with generators $\{e_a\}$ satisfying $[e_a,e_b]=f_{ab}^{\phantom{ab}c}e_c$.
A smooth action of $G$ on $\mfdm$ is a smooth transformation $\Phi:G\times\mfdm \rightarrow \mfdm$ such that denoting $\Phi_g:\mfdm\rightarrow\mfdm$ for every $g\in G$
we have a composition rule compatible with the group structure
$\Phi_g\circ\Phi_h=\Phi_{gh}$.
This induces a pull-back action $\rho$ on the algebra of differential forms by
$\rho_g(\omega):= (\Phi_g^{-1})^{\ast}\omega$ for  $g\in G$ and $\omega\in A$, which we will denote for simplicity as $g\triangleright\omega$. For each $\zeta\in\g$ we use the same symbol for the vector field generating the infinitesimal action of $G$ along $\zeta$ on $\mfdm$. The Lie derivative $L_{\zeta}$ is a degree zero derivation of $A$; denoting by $L_a=L_{e_a}$ the Lie derivatives along generators of $\g$ we have commutation relations $ [L_{e_a},L_{e_b}]= f_{ab}^{\phantom{ab}c}L_{e_c}$ so that $L$ defines a representation of $\g$ on $A$. Thus the algebraic analogue of a $G$ action on $\mfdm$ is a representation of $\g$ on $A$ by derivations; this representation lifts to $\env{\g}$ and the Leibniz rule of $L_{e_a}$ is equivalent to the fact that $e_a$ has primitive coproduct in $\env{\g}$. We will come back on this point when defining covariant actions of Hopf algebras.

We then consider the interior derivative $i_{\zeta}$, defined as the degree $-1$ derivation on $A$ given by contraction along the vector field $\zeta$. In the same way the (infinitesimal) action of $G$ gives a representation of $\g$ (and $\env{\g}$) on $A$, we look now for the algebraic analogue of $i_{\zeta}$.

Out of $\g$ we can construct a super (or $\mathbb{Z}_2$-graded) Lie algebra $\sg=\g\oplus\g$ by adding odd generators $\{\xi_a\}$ that span a second copy of $\g$ as vector space, and putting relations (the brackets are compatible with the degrees)
\begin{equation}
\label{relsg}
[e_a,e_b]=f_{ab}^{\phantom{ab}c}e_c \qquad\qquad [\xi_a,\xi_b]=0 \qquad\qquad [e_a,\xi_b]=f_{ab}^{\phantom{ab}c}\xi_c
\end{equation}
The structure of $\sg$ reflects the usual commutation relations of Lie and interior derivatives; indeed denoting $L_a=L_{e_a}$ and similarly $i_b=i_{e_b}$ it is well known that
\begin{equation}
\label{relil}
[L_a,L_b]=f_{ab}^{\phantom{ab}c}L_c \qquad\qquad [i_a,i_b]=0 \qquad\qquad [L_a,i_b]=f_{ab}^{\phantom{ab}c}i_c
\end{equation}
We can then say that $L_a$ and $i_a$ realize a representation of the super Lie algebra $\sg$ on $A$ as graded derivations; once again this representation lifts to the super enveloping algebra $\env{\sg}$.

To conclude, let us consider also the De Rham differential $d:A^{\bullet}\rightarrow A^{\bullet +1}$
in this algebraic picture. We can add to $\sg$ one more odd generator $d$, obtaining the super Lie algebra
\begin{equation}
\label{gt}
\gt=\sg\oplus \{d\} = \g_{(-1)}\oplus\g_{(0)}\oplus\{d\}_{(1)}
\end{equation}
with relations $(\ref{relsg})$ completed with
\begin{equation}
\label{reld}
[e_a,d]=0 \qquad\qquad [\xi_a,d]=e_a \qquad\qquad [d,d]=0
\end{equation}
The structure induced by $(L,i,d)$ on the algebra of differential forms of a manifold acted by a Lie group may be summarized in the following general definition.
\begin{df}
\label{gdadef}
An algebra $A$ carrying a representation of the super Lie algebra $\gt$ by graded derivations will be called a $\gt$-differential algebra, or $\gt$-da for short.
\end{df}

We pointed out that the fact $(L,i,d)$ act as derivations on differential forms is directly related to the coproduct structure of $\env{\gt}$. The general notion of compatibility between a Hopf algebra $\mathcal{H}$ and the product structure of some algebra $A$ acted by $\mathcal{H}$ is expressed through the definition of covariant actions. A standard reference on Hopf algebras is \cite{maj}, where the omitted definitions and proofs of these introductory sections can be found. We will work with vector spaces, algebras etc.. over the field $\mathbb{C}$.

\begin{df}
\label{covact}
Let $\mathcal{H}$ be a Hopf algebra acting on a unital algebra $A$. The action is said to be covariant if
\begin{equation}
\label{covac}
h\triangleright (ab) := \triangle (h)\triangleright (a\otimes b) = (h_{(1)}\triangleright a)\otimes (h_{(2)}\triangleright b)
\qquad\qquad h\triangleright 1=\epsilon(h)
\end{equation}
When these conditions hold we say that $A$ is a $\mathcal{H}$-module algebra.
\end{df}

\begin{ex}
\label{exadj}
Let $\mathcal{H}$ be a Hopf algebra. A covariant action of $\mathcal{H}$ on itself is given by the left adjoint action
\begin{equation}
\label{adj}
h \triangleright^{ad} g = ad_h(g) = h_{(1)}gS(h_{(2)})
\end{equation}
Note that when $\mathcal{H}=\env{\g}$ for some Lie algebra $\g$ the adjoint action with respect $x\in\g$ equals the bracket with $x$
$$ x\triangleright^{ad}h = ad_x(h) = xh - hx = [x,h] \qquad\qquad x\in\g , \, h\in\env{\g} $$
\end{ex}

\begin{ex}
\label{gdaex}
Let $G$ be a Lie group acting on a manifold $\mfdm$. We already discussed the action of $\g$, $\gt$ and
their enveloping algebras on $A=\Omega^{\bullet}(\mfdm)$, referring to it as a $\g$-da (resp $\gt$-da) structure (see Def(\ref{gdadef})). We now notice that this action is covariant, so the fact that $(L,i,d)$ are (graded) derivations on $A$ is equivalent to the fact that $(e_a,\xi_a,d)$ have primitive coproduct $\triangle(x)=x\otimes 1 + 1\otimes x$. Thus to be a $\gt$-da is equivalent to being a $\env{\gt}$-module algebra.
\end{ex}

We have motivated our interest in the category of (left) Hopf-module algebras, denoted $\hmod$. To study some of its properties in a more efficient language, we present here some basic definitions and facts on braided tensor categories. These ideas are mainly due to Majid; we refer again to his book \cite{maj} for
more details and omitted proofs.

\begin{df}
\label{braidcat}
A braided monoidal (or quasitensor) category $(\mathscr{C}, \otimes, \Psi)$ is a \\ monoidal category $(\mathscr{C},\otimes)$ with a natural equivalence between  the two functors $\otimes, \otimes^{op}:\mathscr{C}\times\mathscr{C}\rightarrow\mathscr{C}$ given by functorial
isomorphisms (called braiding morphisms)
\begin{equation}
\label{braideq}
\Psi_{V,W}: V\otimes W \rightarrow W\otimes V \qquad\qquad \forall \, V,W \, \in \mathscr{C}
\end{equation}
obeying hexagon conditions expressing compatibility of $\Psi$ with the associative structure of $\otimes$ (see for example \cite{maj}(fig $9.4$, pg $430$)). If in addition $\Psi^2=id$ the category $(\mathscr{C}, \otimes, \Psi)$ is said to be a symmetric (or tensor) category.
\end{df}

The relevant example for us is the tensor product of two Hopf-module algebras $A\otimes B$; it is still a Hopf-module algebra, with action defined by
\begin{equation}
\label{tensact}
h\triangleright (a\otimes b) = (h_{(1)}\triangleright a)\otimes (h_{(2)}\triangleright b) \qquad \forall \, a\in A, b\in B, h\in\mathcal{H}
\end{equation}
This means that $\hmod$ is a monoidal category. The algebraic structure of $A\otimes B$ and the presence of a nontrivial braiding operator depend on the quasitriangular structure of $\mathcal{H}$.

\begin{prop}
\label{brcatprop}
If $(\mathcal{H},\mathcal{R})$ is a quasitriangular Hopf algebra the category of left $\mathcal{H}$-module algebras $\hmod$ is a braided monoidal category with braiding morphism
\begin{equation}
\label{brmor}
\Psi_{A,B}(a\otimes b) = (\mathcal{R}^{(2)}\triangleright b)\otimes (\mathcal{R}^{(1)}\triangleright a)
\qquad\qquad \forall \, a\in A, b\in B \; \mbox{and} \; A,B \in \, \hmod
\end{equation}
\end{prop}

Note that when the Hopf algebra is triangular, we may have a non-trivial braiding morphism but it squares to the identity, so that the category is symmetric. If moreover $\mathcal{H}$ is cocommutative, like classical enveloping algebras, $\mathcal{R}=1\otimes 1$ and the braiding morphism is nothing but the flip morphism
$\tau :A\otimes A\rightarrow A\otimes A$ which exchanges the first and second copy of $A$, $\tau(a_1\otimes a_2)=a_2\otimes a_1$.  In this case the ordinary tensor algebra structure of $A\otimes B$, namely
$ (a_1\otimes b_1)\cdot (a_2\otimes b_2) = (a_1a_2) \otimes (b_1b_2)$, is compatible with the action of $\mathcal{H}$. However in the general case, in order to get an algebra structure on $A\otimes B$ acted covariantly by $\mathcal{H}$, we have to take into account the quasitriangular structure; this will be the case for deformed Hopf algebras describing deformed symmetries.

\begin{prop}
\label{braidtens}
If $(\mathcal{H},\mathcal{R})$ is a quasitriangular Hopf algebra and $A, B \in \, \hmod$, the braided tensor product $\mathcal{H}$-module algebra $A\widehat{\otimes}B$ is the vector space $A\otimes B$ endowed with the product
\begin{equation}
\label{braidmult}
(a_1\otimes b_1)\cdot (a_2\otimes b_2) := a_1(\mathcal{R}^{(2)}\triangleright a_2) \otimes (\mathcal{R}^{(1)}\triangleright b_1)b_2
\end{equation}
\end{prop}

The last idea we want to present in this section concerns the notion of commutatitivy; when dealing with a braided category of algebras, it is natural to relate this notion to the braiding morphism of the category. Indeed the commutatitivy of an algebra $A$ may be expressed as the commutativity of the multiplication $m_A:A\otimes A$ with the flip morphism $\tau$; when we are interested in a specific category, in our case $\hmod$, it is natural to ask that both the maps are morphism in the category. The multiplication map $m_A$ is a morphism in $\hmod$ excatly by definition of covariant action, while for $\mathcal{H}$ quasitriangular we know that $\tau$ is no longer a morphism in $\hmod$, but its natural analogue is the braiding morphism $\Psi$. This motivates the following definition.

\begin{df}
\label{braidcomm}
In the category $\hmod$ an algebra $A$ is said to be braided commutative if its multiplication map 
$m:A\otimes A\rightarrow A$ commutes with the braiding morphism $\Psi_{A,A}$:
\begin{equation}
\label{brgr}
m\circ\Psi_{A,A} = m \qquad\qquad \Longleftrightarrow
\qquad\qquad a\cdot b = (\mathcal{R}^{(2)}\triangleright b)\cdot (\mathcal{R}^{(1)}\triangleright a)
\end{equation}
\end{df}
Thus the property to be commutative now depends on the Hopf algebra which acts; it could happen that an algebra is acted covariantly by two different Hopf algebras and it is braided commutative with respect the first one but not with respect the second one. 

\subsection{Deformation of symmetries by Drinfel'd twists}

Using the language of Def (\ref{gdadef}) we will consider a symmetry acting on a graded algebra $A$ as being expressed by a $\gt$-da structure on $A$. By deformation of a symmetry we mean a deformation of the Lie algebra $\gt$ or a deformation of the Hopf algebra $\env{\gt}$. To the first case belong quantum Lie algebras, while the second case refers to quantum enveloping algebras.

In both the approaches, and depending on the particular deformation considered, a general strategy is to relate the deformation of $\gt$ or $\env{\gt}$ to a deformation of the product in every $\gt$-da $A$, and vice versa. When such a link between symmetries (i.e. Hopf or Lie algebras), spaces (i.e. $\gt$-da) and deformations is present, we will speak of covariant deformations or induced star products.

We give a detailed presentation of this ideas by picking up a particular class of deformations, the ones generated by Drinfel'd twists in Hopf algebras \cite{dr1,dr2}; we choose to work with Drinfel'd twists for several reasons. They provide the most natural setting to describe and study symmetries of a large class of nc geometries, like toric isospectral deformations, Moyal planes or nc toric varieties, they allow for quite explicit computations and moreover they often are the only class of deformations up to isomorphism, as we breifly discuss at the end of the section. 

Thus the following exposition will be focused on this specific, even if quite general, class of deformations. However we feel that the general strategy to study nc actions and define nc equivariant cohomology is actually independent from the specific deformation choosen, thus part of what we are going to present could in principle be applied to different class of deformations; we will say more on this in section $2.5$. 

We start with the definition and basic properties of Drinfel'd twists. For omitted proofs and a more detailed exposition we remand to the original works of Drinfel'd \cite{dr1}\cite{dr2} or to \cite{maj}.

\begin{df}
\label{twisteldf}
Let $\mathcal{H}$ be an Hopf algebra. An element $\chi=\chi^{(1)}\otimes\chi^{(2)} \in\mathcal{H}\otimes\mathcal{H}$ is called a twist element for $\mathcal{H}$ if it satisfies the following properties:
\begin{enumerate}
\item $\chi$ is invertible
\item $(1\otimes \chi)(id\otimes\bigtriangleup)\chi = (\chi\otimes 1)(\bigtriangleup\otimes id)\chi$ \hspace{1cm} (cocycle condition)
\item $(id\otimes\epsilon)\chi=(\epsilon\otimes id)=1$ \hspace{1cm} (counitality)
\end{enumerate}
\end{df}

\begin{thm}
\label{twistthm}
A twist element $\chi=\chi^{(1)}\otimes\chi^{(2)}\in\mathcal{H}\otimes\mathcal{H}$ defines a twisted Hopf algebra structure $\mathcal{H}^{\chi}=(\mathcal{H},\cdot,\triangle^{\chi},S^{\chi},\epsilon)$ with same multiplication and counit and new coproduct and antipode given by
\begin{equation}
\triangle^{\chi}(h)=\chi\triangle(h)\chi^{-1} \quad , \quad S^{\chi}(h)=US(h)U^{-1} \quad \mbox{ with }
\quad U=\chi^{(1)}S\chi^{(2)}
\end{equation}
When applied to quasitriangluar Hopf algebras $(\mathcal{H},\mathcal{R})$ the twist deforms the quasitriangular structure to $\mathcal{R}^{\chi}=\chi_{21}\mathcal{R}\chi^{-1}$ ($\chi_{21}=\chi^{(2)}\otimes\chi^{(1)}$).
\end{thm}

We point out that the cocycle condition on $\chi$ is a sufficient condition to preserve the coassociativity of the coproduct. A more general theory of twists where this requirement is dropped out is well defined in the category of quasi-Hopf algebras \cite{dr1}\cite{dr2}. The theory of Drinfel'd twists easily extends to super (or $\mathbb{Z}_2$ graded) Hopf algebras; this will be relevant for our purposes, since we are interested in deformations of $\env{\gt}$. 

\begin{thm}
\label{covdef}
If $A$ is a left $\mathcal{H}$-module algebra and $\chi$ a Drinfeld twist for $\mathcal{H}$, the deformed product
\begin{equation}
\label{defprod}
a\cdot_{\chi} b:=\cdot\left(
\chi^{-1}\triangleright(a\otimes b)\right) \qquad \qquad  \forall \,
a,b\in A
\end{equation}
makes $A_{\chi}=(A,\cdot_{\chi})$ into a left $\mathcal{H}^{\chi}$-module algebra with respect to the same action. 
\end{thm}

Thus a Drinfel'd twist in $\mathcal{H}$ generates a deformation of the algebra structure of every $\mathcal{H}$-module algebra; by interpreting deformed module-algebras as nc spaces, we may think of $\mathcal{H}$ (or its Drinfel'd twist element) as the source of noncomutativity. 

There is also a dual notion of Drinfel'd twists \cite{maj}, where the multiplication of $\mathcal{H}$ is deformed while the coproduct is unchanged; in this case the induced deformation involves Hopf-comodule algebras. Since we prefer to work with actions of $\gt$ and deformation of its enveloping algebra we will use the Drinfel'd twist of Thm \ref{twistthm}, but everything could be restated in term of coactions of the Hopf algebra of (representable) functions over the group $G$ (which is dual to $\env{\g}$) and its dual Drinfel'd twist.

So we focus on Drinfe'ld twists of enveloping algebras $\env{\g}$. In order to have more explicit computations we restrict to the case of semisimple Lie algebras, so that we have at our disposal a Cartan decomposition of $\g$ with an abelian Cartan subalgebra $\h$. Moreover we use twist elements $\chi$ contained in $\env{\h}\otimes\env{\h}\subset\env{\g}\otimes\env{\g}$; we refer to this choice as the class of abelian Drinfeld twists, in the sense that $[\chi,\chi]=0$. A general theory for Drinfeld twist deformations of enveloping algebras with non abelian twist elements could lead to very interesting results, and deserves a detailed study in the future.

After these assumptions, let us fix the notations. Given a semisimple Lie algebra $\g$ we fix a Cartan decomposition
$$\{H_i,E_r\}  \qquad\qquad i=1,\ldots ,n , \quad r=(r_1,\ldots ,r_n) \in\mathbb{Z}^n$$
where $n$ is the rank of $\g$, $H_i$ are the generators of the Cartan subalgebra $\h\subset\g$ and $E_r$ are the roots element labelled by the $n$-dimensional root vectors $r$. In this decomposition the structure constants are written as follows:
\begin{equation}
\label{grel}
\begin{array}{rclcrcl}
           [H_i, H_j]   & = & 0              & &  [H_i,E_r] & = & r_i E_r \\
\phantom{3}[E_{-r},E_r] & = & \sum_i r_i H_i & &  [E_r,E_s] & = & N_{r,s}E_{r+s}
\end{array}
\end{equation}
The explicit expression of $N_{r,s}$ is not needed in what follows, but it worths saying that it vanishes if $r+s$ is not a root vector.

Now we choose a twist element $\chi$ which depends on Cartan generators $H_i$. Since we use the Drinfel'd twist as a source of 'quantization' or deformation, we want it to depend on some real parameter(s) $\theta$ and recover the classical enveloping algebra for $\theta\rightarrow 0$. Thus we are actually making a Drinfeld twist in the formal quantum enveloping algebra $\env{\g}_{[[\theta]]}$.
We will make use of the following twist element, which first appeared in \cite{res}:
\begin{equation}
\label{twistel}
\chi = \mbox{exp }\{-\frac{i}{2}\, \theta^{kl}H_k\otimes H_l\} \qquad\qquad\qquad\qquad \chi\in(\env{\h}\otimes\env{\h})_{[[\theta]]}
\end{equation}
with $\theta$ a $p\times p$ real antisymmetric matrix, $p\leq n$ (i.e. we do not need to use the whole $\h$ to generate the twist).

Using relations (\ref{grel}) and the expressions in Thm \ref{twistthm} for the twisted coproduct and antipode, we can describe explicitly the Hopf algebra structure of $\envt{\g}_{[[\theta]]}$.

\begin{prop}
\label{twcopr}
Let $\chi$ be the twist element in (\ref{twistel}). The twisted coproduct $\triangle^{\chi}$ of $\envt{\g}_{[[\theta]]}$ on the basis $\{H_i,E_r\}$ of $\g$ reads
\begin{align}
\label{copH}\triangle^{\chi}(H_i) & = \triangle (H_i) = H_i\otimes 1 + 1\otimes H_i\\
\label{copE}\triangle^{\chi}(E_r) & = E_r\otimes\lambda_r^{-1} + \lambda_r\otimes E_r
\end{align}
where
\begin{equation}
\label{lambda}
\lambda_r=\mbox{exp }\{\frac{i}{2}\,\theta^{kl}r_kH_l\}
\end{equation}
are group-like elements (one for each root $r$) with untwisted coproduct $\triangle^{\chi}(\lambda_r)=\triangle (\lambda_r) = \lambda_r\otimes\lambda_r$.
\end{prop}
\proof{It is clear that $\forall X \in \mathcal{H}$ whenever $[H_i,X]=0$ the coproduct $\triangle(X)$ is not deformed. Thus $(\ref{copH})$ follows easily; for $(\ref{copE})$ we compute
$$ \mbox{exp }\{-\frac{i}{2}\, \theta^{\mu\nu}H_{\mu}\otimes H_{\nu}\} (E_r\otimes 1 + 1\otimes E_r)
   \, \mbox{exp }\{ \frac{i}{2}\, \theta^{\mu\nu}H_{\mu}\otimes H_{\nu}\} $$
at various order in $\theta$, using
$$e^{tA}Be^{-tA}=\sum_{n=0}^{\infty} \frac{t^n}{n!}[A,[A,\ldots [A,B]]]$$
At the first order we have
\begin{equation*}
\begin{split}
-\frac{i}{2}\theta^{\mu\nu}[H_{\mu}\otimes H_{\nu}, E_r\otimes 1 + 1\otimes E_r] & =
-\frac{i}{2}\theta^{\mu\nu}\left( [H_{\mu},E_r]\otimes H_{\nu} + H_{\mu}\otimes [H_{\nu},E_r] \right) \\
 & = -\frac{i}{2}\theta^{\mu\nu} \left( E_r\otimes r_{\mu}H_{\nu} + r_{\nu}H_{\mu}\otimes E_r \right)
\end{split}
\end{equation*}
So the second order is
\begin{equation*}
\begin{split}
(\frac{i}{2})^2\theta^{\mu\nu}\theta^{\rho\sigma} & [H_{\mu}\otimes H_{\nu}, E_r\otimes r_{\rho}H_{\sigma} + r_{\sigma}H_{\rho}\otimes E_r] = \\
 & = (\frac{i}{2})^2\theta^{\mu\nu}\theta^{\rho\sigma} \left( [H_{\mu},E_r]\otimes r_{\rho}H_{\nu}H_{\sigma} +
    r_{\sigma}H_{\mu}H_{\nu}\otimes [H_{\nu},E_r] \right) = \\
 & = (\frac{i}{2})^2\theta^{\mu\nu}\theta^{\rho\sigma} \left( E_r\otimes r_{\mu}r_{\rho}H_{\nu}H_{\sigma} +
    r_{\sigma}r_{\nu}H_{\mu}H_{\rho}\otimes E_r \right)
\end{split}
\end{equation*}
and carrying on with higher orders the series gives (\ref{copE}).}

\begin{prop}
\label{twantip}
Let $\chi$ be the twist element in (\ref{twistel}). The element $U=\chi^{(1)}S\chi^{(2)}$ reduces to the identity so that the twisted antipode $S^{\chi}(h)=US(h)U^{-1}$ equals the untwisted one.
\end{prop}
\proof{We compute $U$ at various order in $\theta$. The order zero is trivially the identity; the first order is
\begin{equation*}
- \frac{i}{2}\theta^{\mu\nu}H_{\mu}S(H_{\nu}) = \frac{i}{2}\theta^{\mu\nu}H_{\mu}H_{\nu}
\end{equation*}
and so it vanishes by antisymmetry of $\theta^{\mu\nu}$. The same happens to the second order
\begin{equation*}
(\frac{i}{2})^2\theta^{\mu\nu}\theta^{\rho\sigma}H_{\mu}H_{\rho}S(H_{\nu}H_{\sigma}) =
(\frac{i}{2})^2\theta^{\mu\nu}\theta^{\rho\sigma}H_{\mu}H_{\rho}H_{\sigma}H_{\nu} = 0
\end{equation*}
and all higher orders are zero for the same reason.}

Finally, the twisted quasitriangular structure (we start with $\mathcal{R}=1\otimes 1$ in $\env{\g}$) is
\begin{equation}
\mathcal{R}^{\chi}=\chi_{21}\mathcal{R}\chi^{-1}=\chi^{-1}(1\otimes 1)\chi^{-1}=\chi^{-2}
\end{equation}
so the twisted enveloping algebra is triangular but no more cocommutative. This completes the explicit computation of the Hopf algebra structure of $\envt{\g}_{[[\theta]]}$. \\

We end this section with a brief discussion on the relation between Drinfel'd twists and other deformations of enveloping algebras; we refer to \cite{kas}\cite{ss} for a detailed treatment and the proofs. The theory of algebras and coalgebras deformations, and related cohomologies, is well defined  in the setting of formal power series; the results we quickly present here are mainly due to Gerstenhaber, Schack, Shnider and Drinfel'd.

To introduce quantum enveloping algebras several approaches are possible: a first possibility is to consider deformations $\g_{\theta}$ of the Lie algebra structure of $\g$, basically by defining structure constants on $\mathbb{C}_{[[\theta]]}$, so that $(\env{\g_{\theta}},\cdot_{\theta},\triangle_{\theta},\mathcal{R}_{\theta})$ is the associated quantum enveloping algebra defined using the $\theta$-deformed brackets in $\g_{\theta}$. However a classical result in deformation theory, due to Gerstenhaber, states that if an algebra $A$ has a vanishing second Hochschild cohomology group $H^2(A,A)=0$, then any deformation $A^{\prime}$ is isomorphic to the $\theta$-adic completion of the undeformed algebra, i.e., $A^{\prime}\simeq A_{[[\theta]]}$; these algebras are called rigid. For example for semisimple Lie algebras rigidity is implied by the second Whitehead lemma, and so they only admit trivial deformations.

When $\g$ is semisimple a standard deformation of its enveloping algebra is provided by the Drinfel'd-Jimbo quantum enveloping algebra $\mathfrak{U}_{\theta}(\g)$, defined as the topological algebra over $\mathbb{C}_{[[\theta]]}$ generated by Cartan and roots element $\{H_i,X_i,Y_i\}$ subjects to relations ($a_{ij}$ is the Cartan matrix and $D=(d_1\ldots d_n)$ the diagonal matrix of root length)
\begin{align}
[H_i,H_j] & = 0         & [X_i,Y_j]& = \delta_{ij} \frac{\sinh(\theta d_iH_i/2)}{\sinh(\theta d_i/2)} \\
[H_i,X_j] & = a_{ij}X_j & [H_i,Y_j]& = -a_{ij}Y_j
\end{align}
plus the $\theta$-quantized version of Serre relations between $X_iX_j$ and $Y_iY_j$ for $i\neq j$.
Now, the rigidity of $\g$ assures that there is an isomorphism of topological algebras
$$\alpha:\mathfrak{U}_{\theta}(\g)\rightarrow\env{\g}_{[[\theta]]}$$
which transfers the Hopf algebra structure $\triangle_{\theta},\epsilon_{\theta},S_{\theta}$ of $\mathfrak{U}_{\theta}(\g)$
to $\env{\g}_{[[\theta]]}$ by
\begin{equation}
\label{prhopf}
\triangle^{\prime} = (\alpha\otimes\alpha)\circ\triangle_{\theta}\circ\alpha^{-1} \quad , \quad
\epsilon^{\prime} = \epsilon_{\theta}\circ\alpha^{-1} \quad , \quad
S^{\prime} = \alpha\circ S_{\theta}\circ\alpha^{-1}
\end{equation}
so that $\alpha$ becomes an isomorphism of Hopf algebras from $\mathfrak{U}_{\theta}(\g)$ to $\env{\g}_{[[\theta]]}$ (with the primed Hopf algebra
structure of (\ref{prhopf})). Now, again for rigidity reasons the two coproducts $\triangle$ and $\triangle^{\prime}$ in $\env{\g}_{[[\theta]]}$ must be related by an inner automorphism: there should exist an invertible element $\chi\in(\env{\g}\otimes\env{\g})_{[[\theta]]}$ such that $\triangle^{\prime}(h)=\chi\triangle(h)\chi^{-1}$. This $\chi$ quite often does not satisfy any cocycle condition, so it defines a generalized Drinfel'd twist and $\envt{\g}_{[[\theta]]}$ is a quasi-Hopf algebra with a nontrivial coassociator $\Phi$ encoding basically all the information about the Drinfel'd-Jimbo deformation.

So, at least for rigid Lie algebras, there is only one class of deformations modulo isomorphism. We can equivalently consider either deformations involving Lie algebra generators and their relations, as in the spirit of $\mathfrak{U}_q(\g)$, or rather (generalized) Drinfel'd twists of $\env{\g}_{[[\theta]]}$ in which the algebra structure is undeformed and the whole deformation is contained in the coproduct (plus eventually a non trivial coassociator).

\subsection{Toric isospectral deformations from Drinfel'd twists}

In the previous section we fixed the class of Drinfel'd twists $\chi$ we are interested in (\ref{twistel}), and noted that they are generated by elements in the Cartan subalgebra $\h$ of a semisimple Lie algebra $\g$. Then we showed that as a consequence of the twist every $\env{\g}$-module algebra deforms its product in order to preserve the covariance of the action. Following this strategy, it is clear we can induce a nc deformation in the algebra of functions (or differential forms) of every manifold acted upon by some group of rank $\geq 2$.\footnote{With only one toric generator the twist element \ref{twistel} is necessarily trivial, i.e. a coboundary. See \cite{maj}(Prop $2.3.5$) for a cohomological classification of Drinfel'd twists.}

This is the setting of toric isospectral deformations \cite{cl}\cite{cdv}. One starts with a compact Riemannian spin manifold $\mfdm$ whose isometry
\footnote{In the construction of the deformed spectral triple this property is crucial since it assures the invariance of the Dirac operator. This fact however does not concern the deformation of the algebra $C^{\infty}(\mfdm)$, and so we can relax this request in the Drinfeld twist approach. Nevertheless note that the action of a compact Lie group $G$ on a Riemannian manifold $(\mfdm,g)$ can always turned into an isometry by averaging the metric $g$ with respect to the action of the group.}
group has rank at least $2$, and uses the action of the compact torus $\mathbb{T}^n$ ($n\geq 2$) to construct a nc spectral triple $(C^{\infty}(\mfdm_{\theta}),L^2(\mfdm,S),D)$ by deforming the classical one; the name 'isospectral' refers to the fact that in the nc spectral triple the algebra of functions and its representation on $L^2(\mfdm,S)$ are deformed, but not the Dirac operator $D$ (and so its spectrum) which is still the classical one due to its invariance under the action.

We quickly review the construction of toric isospectral deformations, and then we show that the same algebra deformation can be obained by a Dinfel'd twist in the enveloping algebra of the torus. We do not discuss the full spectral triple of toric isospectral deformations, since our interest is contained in the category of Hopf-module algebras; however it is implicit that when we say we can interpret deformed Hopf-module algebras as nc spaces we have to describe the whole spectral triple to give a full meaning to the name 'induced nc geometry'.

Under the hypothesis of compactness of $\mfdm$ we can decompose the algebra of smooth functions $C^{\infty}(\mfdm)=\bigoplus_{r\in (\mathbb{Z}^n)^{\ast}}C^{\infty}_r(\mfdm)$ in spectral subspaces labelled by weights $r$ of the torus action, such that every $f_r\in C^{\infty}_r(\mathcal{M})$ is an eigenfunction. Representing elements of $\mathbb{T}^n$ as $e^{2\pi it}$ with $t\in\mathbb{Z}^p$, the action $\sigma$ on an eigenfunction $f_r$ is given by a phase factor depending on $r$:
\begin{equation}
\sigma_t(f_r) = e^{2\pi i t\cdot r}f_r \qquad \qquad \qquad t\in\mathbb{Z}^n \, , \; r \in (\mathbb{Z}^n)^{\ast}
\end{equation}
Taking a real $n\times n$ skew-symmetric matrix $\theta$ we can define a deformed product between eigenfunctions
\begin{equation}
  \label{isosp}
  f_r\times_{\theta}g_s:=\mbox{exp }[\frac{i}{2}\,\theta^{kl}r_ks_l] f_rg_s
\end{equation}
and by linearity extend it on the whole of $C^{\infty}(\mathcal{M})$. We will call
\begin{equation}
\label{ncalg}
\ncf :=\left(C^{\infty}(\mathcal{M}), \times_{\theta}\right)
\end{equation}
the algebra of functions of the nc manifold $\mathcal{M}_{\theta}$. Clearly, $\mathbb{T}^n$-invariant functions form a commutative ideal in the nc algebra $C^{\infty}(\mfdm_{\theta})$.

The deformed product (\ref{isosp}) is a sort of Moyal product, with the action of $\mathbb{R}^n$ replaced by the torus $\mathbb{T}^n$, i.e. considering periodic actions of $\mathbb{R}^n$. The idea to use actions (of $\mathbb{R}^n$) to produce strict deformation quantizations indeed appeared firstly in \cite{rie}.

We now express the previous deformation in the language of Drinfel'd twists. Since we supposed the compact Lie group $G$ acting on $\mfdm$ to have rank $n\geq 2$, we can use its Cartan generators $H_i\in\h\subset\g$ $(i=1,\ldots n)$ and the real skewsymmetric matrix $\theta$ to define a twist element $\chi\in (\env{\g}\otimes\env{\g})_{[[\theta]]}$ (the same of (\ref{twistel}))
\begin{equation*}
\chi = \mbox{exp }\{-\frac{i}{2}\, \theta^{kl}H_k\otimes H_l\}
\end{equation*}
We already computed the twisted Hopf algebra structure of $\envt{\g}_{[[\theta]]}$ in section 1.2; now following Thm \ref{covdef} we describe the deformed product induced on the $\env{\g}$-module algebra $A=\Omega(\mfdm)$. As we did for functions, we decompose $A=\oplus_r\mathcal{A}_r$ into spectral subspaces labelled by characters of the toric subgroup of $G$ so that $H_k\triangleright\omega_r=r_k\omega_r$. On the spectral subspaces the induced deformed product is easily computed.

\begin{prop}
On spectral elements $\omega_r\in A_r$ and $\omega_s\in A_s$ the product induced from the Drinfeld twist of $\env{\g}$ reads
\begin{equation}
\label{ncdf}
\omega_r\wedge_{\theta}\omega := \chi^{-1}\triangleright (\omega_r\otimes\omega_s)=
\mbox{exp }\{\frac{1}{2}\theta^{\mu\nu}r_{\mu}s_{\nu}\} \omega_r\wedge\omega_s
\end{equation}
\end{prop}
\proof{The result follows from a direct computation, using the explicit expression of $\chi$ and
$$\theta^{\mu\nu}(H_{\mu}\otimes H_{\nu}) \triangleright (\omega_r\otimes\omega_s) = \theta^{\mu\nu}r_{\mu}s_{\nu}(\omega_{\mu}\otimes\omega_{\nu})$$
which use the spectral property of $\omega_r$ and $\omega_s$.}

We extend this product from spectral elements to the whole algebra $A$ by linearity.
\begin{df}
\label{ncdfdef}
The nc algebra $A_{\chi}=(A, \, \wedge_{\theta})$ with product $\wedge_{\theta}$ defined in (\ref{ncdf}) is called the algebra of nc differential forms of the nc space $\mfdm_{\theta}$.
\end{df}

The degree zero part of $A_{\chi}$ is the algebra $C^{\infty}(\mfdm_{\theta})$ of (\ref{ncalg}). This shows it is possible to recover toric isospectral algebra deformations by Drinfel'd twists.

We deformed the graded commutative wedge product $\wedge$ to obtain a nc product $\wedge_{\theta}$.
Recalling Def \ref{braidcomm} (and the natural generalization to graded-commutative algebras) a natural question is then if $\wedge_{\theta}$ is braided graded-commutative. 

\begin{lemma}
\label{l1}
Let $A$ be a graded commutative algebra in $\hmod$ and $\chi$ a twist element of the form (\ref{twistel}).
Then
\begin{equation}
a_1\cdot_{\chi}a_2 := \cdot (\chi^{-1}\triangleright (a_1\otimes a_2)) = (-1)^{|a_1||a_2|} \cdot (\chi\triangleright (a_2\otimes a_1))
\end{equation}
\end{lemma}
\noindent\underline{\textsf{Proof}}: By direct computation, starting from the rhs:
\begin{equation*}
\begin{split}
(-1)^{|a_1||a_2|} & \left( \sum_n (-\frac{\theta^{\alpha\beta}}{2})^n \frac{1}{n!} \, (H_{\alpha}^na_2)\cdot (H_{\beta}^na_1) \right) =
 \sum_n (-\frac{\theta^{\alpha\beta}}{2})^n \frac{1}{n!} \, (H_{\beta}^na_1)\cdot (H_{\alpha}^na_2)  = \\
 & = \sum_n ( \frac{\theta^{\beta\alpha}}{2})^n \frac{1}{n!} \, (H_{\beta}^na_1)\cdot (H_{\alpha}^na_2)  =
 \sum_n ( \frac{\theta^{\alpha\beta}}{2})^n \frac{1}{n!} \, (H_{\alpha}^na_1)\cdot (H_{\beta}^na_2) = \\
 & = \cdot \, (\chi^{-1} \triangleright (a_1\otimes a_2)) = a_1\cdot_{\chi}a_2 \qquad\qquad\qquad _\blacksquare
\end{split}
\end{equation*}

\begin{prop}
\label{propgrc}
Let $A_{\chi}$ be the algebra of nc differential forms deformed by the usual Drinfeld twist (\ref{twistel}) in $\envt{\g}_{[[\theta]]}$. Then $A_{\chi}$ is braided graded-commutative (see Def \ref{braidcomm}).
\end{prop}

\noindent\underline{\textsf{Proof}}: The quasitriangular structure of $\envt{\gt}$ is $\mathcal{R}^{\chi}=\chi^{-2}$. We compute the rhs of (\ref{brgr}) with $\omega\in A_{\chi}^n$ and $\nu\in A_{\chi}^k$, and make use of the previous Lemma:
\begin{equation*}
\begin{split}
(-1)^{kn} ((\mathcal{R}^{\chi})^{(2)} & \triangleright \nu) \wedge_{\theta} ((\mathcal{R}^{\chi})^{(1)}\triangleright \omega) = \\
& = (-1)^{kn} \wedge \left( (\mathcal{R}^{\chi})^{(2)}\otimes (\mathcal{R}^{\chi})^{(1)} \cdot \chi^{-1} \triangleright (\nu\otimes\omega) \right) = \\
 & = (-1)^{kn} \wedge \left(\chi^2\cdot\chi^{-1} \triangleright (\nu\otimes\omega) \right) =
(-1)^{kn} \wedge (\chi\triangleright (\nu\otimes\omega) ) =  \\
 & = \wedge \, (\chi^{-1}\triangleright (\omega\otimes\nu) ) = \omega\wedge_{\theta}\nu \qquad\qquad\qquad _\blacksquare
\end{split}
\end{equation*}

We presented the result having in mind the deformed product in the algebra of differential forms, but it should be clear that the same conclusion applies to every graded-commutative algebra $A$ deformed using a Drinfeld twist of the form (\ref{twistel}) starting from a cocommutative Hopf algebra; in all these cases the deformed product in $A_{\chi}$ turns out to be braided graded-commutative. 

We can summarize this result by saying that in this kind of induced nc deformations all the information regarding the noncommutativity may be encoded in the braiding of the category $\hmod$. Thus toric isospectral deformations and similar nc spaces may be in some sense thought as commutative spaces, only in a braided category; it worths saying that this philosophy to turn properties of objects by shifting the point of view and changing the category in which to look at them is mainly due to Majid, and named by him transmutation \cite{maj}. \\

We conclude this section by showing explicitly how to deform symmetries in order to have a  
covariant action (i.e. a $\gt$-da structure) on a nc algebra of the type (\ref{ncdf}); this means we are now starting with a nc space where the deformation comes from some Dinfel'd twist $\chi$, thus whose nc algebra will be denoted $A_{\chi}$, and we want to accordingly deform every symmetry acting on $A$. The idea is simple and come directly from Thm \ref{covdef}. Our undeformed symmetry was some $\gt$-da struture on $A$; in order to act covariantly on $A_{\chi}$ we know we have then to make a Drinfel'd twist by $\chi$ on $\env{\gt}$. The Lie and interior derivative along generators which commute with $\chi$ will remain undeformed, while the others will satisfy a twisted Leibniz rule due to the deformed coproduct. We are going to show the explicit formulas, which easily follows from Prop \ref{twcopr}. Before doing that, a small remark: we said we have to twist $\env{\gt}$ with the same $\chi$ which deforms the space. In general $\chi$ does not need to belong to $\env{\gt}\otimes\env{\gt}$; we are actually twisting the enveloping algebra of $\tilde{\mathfrak{g}}^{\prime}=\gt\rtimes\kkt$ where $\kkt$ is the symmetry to which $\chi$ belongs and the structure of semidirect product depends on the action of $\kk$ on $\g$. For example if $\g$ and $\kk$ commute we have $\envt{\gt}\cong\env{\gt}$ and the symmetry is undeformed. For simplicity we will directly assume that $\kk\subset\g$; when this is not the case, we need just to replace everywhere $\gt$ with $\tilde{\mathfrak{g}}^{\prime}$.

\begin{df}
\label{defsymm}
A deformed symmetry on a nc algebra $A_{\chi}$ is a twisted $\gt$-da structure, i.e. a covariant action of $\envt{\gt}$. The generators $\{\xi_a,e_a,d\}$ of $\gt$ represent respectively interior derivative, Lie derivative and de Rham differential.
\end{df}

We have already computed the twisted Hopf structure of $\env{\g}$; it remains to describe the action of the twist on the odd part. Following our usual notation, with $\xi_i$ we mean generators corresponding to Cartan-type indexes while  $\xi_r$ denotes root-type indexes.

\begin{prop}
\label{twodd}
The twisted coproduct on odd generators $\{\xi_i,\xi_r\}$ reads
\begin{align}
\label{copxi} \triangle^{\chi}(\xi_i) & = \triangle (\xi_i) = \xi_i\otimes 1 + 1\otimes \xi_i\\
\label{copxir}\triangle^{\chi}(\xi_r)   & = \xi_r\otimes\lambda_r^{-1} + \lambda_r\otimes \xi_r
\end{align}
The twisted antipode $S^{\chi}(\xi_a)$ is equal to the untwisted one, both for Cartan and root generators.
\end{prop}
\proof{For the coproduct part, the proof is just like in Prop \ref{twcopr}; one computes explicitly the definition of $\triangle^{\chi}(\xi_a)$ and use commutation relations between $\xi_a$ and $H_i$. For the antipode, we already showed in Prop \ref{twantip}  that the element $U$ entering in the definition of $S^{\chi}$ for this class of Drinfeld twists is the identity, and so the antipode is undeformed regardless of whether it is computed on even or odd generators.}

We are now ready to interpret the above results in terms of deformed Lie and interior derivatives.

\begin{prop}
\label{ltw}
The Lie derivative $L_{e_a}$ acts classically on single generators of $A_{\chi}$; on the product of two generators $\omega ,\eta \in A_{\chi}$ it satisfies a deformed Leibniz rule:
\begin{align}
\label{li}
L_{H_i}(\omega\wedge_{\theta}\eta) & = (L_{H_i}\omega)\wedge_{\theta}\eta + \omega\wedge_{\theta}(L_{H_i}\eta) \\
\label{le}
L_{E_r}(\omega\wedge_{\theta}\eta) & = (L_{E_r}\omega)\wedge_{\theta}(\lambda_r^{-1}\triangleright\eta) +
(\lambda_r\triangleright\omega)\wedge_{\theta}(L_{E_r}\eta)
\end{align}
For this reason we call $L_{E_r}$ a twisted degree $0$ derivation of the algebra $A_{\chi}$.
\end{prop}
\proof{By definition $L_{e_a}(\omega)=e_a\triangleright\omega$; the claimed formulas are just a restatement of the $\envt{\gt}$-module structure of $A_{\chi}$ taking into account the twisted coproduct (\ref{copH})(\ref{copE}). Note that $\lambda_r\triangleright\omega$ involves only Lie derivatives along Cartan generators.}

\begin{prop}
\label{itw}
The interior derivative $i_a=i_{\xi_a}$ acts undeformed on single generators of $A_{\chi}$;
on products of differential forms it satisfies a deformed graded Leibniz rule:
\begin{align}
\label{ii}
i_{\xi_i}(\omega\wedge_{\theta}\eta) & = (i_{\xi_i}\omega)\wedge_{\theta}\eta + (-1)^{|\omega|}\omega\wedge_{\theta}(i_{\xi_i}\eta) \\
\label{ie}
i_{\xi_r}(\omega\wedge_{\theta}\eta) & = (i_{\xi_r}\omega)\wedge_{\theta}(\lambda_r^{-1}\triangleright\eta) + (-1)^{|\omega|}
(\lambda_r\triangleright\omega)\wedge_{\theta}(i_{\xi_r}\eta)
\end{align}
For this reason $i_{\xi_r}$ is called a twisted derivation of degree $-1$ of the algebra $A_{\chi}$.
\end{prop}
\proof{By definition $i_{\xi_r}(\omega)=\xi_r\triangleright\omega$. The proof is the same of Prop \ref{ltw}, now using the twisted coproduct of odd generators presented in Prop \ref{twodd}.}

The differential $d$ is completely undeformed, since it commutes with the generators of the twist $\chi$ and thus $\triangle^{\chi}(d)=\triangle(d)$. One can also check directly from the definition of $\wedge_{\theta}$ that $d$ satisfies the classical Leibniz rule.

Note that since the Drinfel'd twist in $\env{\gt}$ does not change the Lie brackets in $\gt$, i.e. the Lie algebra structure of $\gt$ is undeformed, the twisted derivations $(L,i,d)$ still obey the classical commutation relations (\ref{relil}).

\begin{ex}
To clarify the relation between the generators of the twist $\chi$ and the symmetry eventually deformed, we consider rotations on the Moyal plane. Similarly to toric isospectral deformations, the Moyal plane $\mathbb{R}^{2n}_{\Theta}$ may be described by a nc algebra deformed by a Drinfel'd twist of the form (\ref{twistel}) but where now the toric generators $H_i$ are replaced by momenta $P_i$
\cite{rie}. Deformed rotations on $\mathbb{R}^{2n}_{\Theta}$, accordingly to Def \ref{defsymm}, are described by a twist of the enveloping algebra $\env{\widetilde{\mathfrak{so}(2n)}}$; since the translations $P_i$ which generates the twist do not belong to the symmetry $\mathfrak{so}(2n)$, this is a situation where we must consider the enveloping algebra of the semidirect product $\mathfrak{so}(2n)\rtimes\mathbb{R}^{2n}$, i.e. of the euclidean group $\mathfrak{e}_{2n}$. Thus denoting by $M_{\mu\nu}$ the generators of $\mathfrak{so}(2n)$ from $[M_{\mu\nu},P_{a}]=g_{\mu a}P_{\nu}-g_{\nu a}P_{\mu}$
we get the twisted coproduct
\begin{equation}
\label{twistrot}
\triangle^{\chi}(M_{\mu\nu})=\triangle(M_{\mu\nu}) + \frac{i\Theta^{ab}}{2}[(\delta_{\mu a}P_{\nu}-\delta_{\nu a}P_{\mu})\otimes P_b + P_a\otimes (\delta_{\mu b}P_{\nu} - \delta_{\nu b}P_{\mu})]
\end{equation}
This means that Lie and interior derivatives along generators of rotations, when acting on a product, satisfy a deformed Leibniz rule which contains extra-terms involving translations.
\end{ex}


\section{Models for noncommutative equivariant cohomology}


The subject of this section is to introduce algebraic models for the equivariant cohomology of nc spaces acted by deformed symmetries. We will do it showing how to recover Weil and Cartan models in the deformed case.

In the first subsection we review some classical notions of equivariant cohomoloy, underliying the role played by the Weil algebra. In the second subsection we describe the nc Weil algebra introduced by Alekseev and Meinrenken \cite{am1}\cite{am2} and their Weil and Catan models for what they call nc equivariant cohomology. In the third subsection we show how to adapt these constructions to the class of nc spaces we described so far, arriving to the definition of a twisted nc equivariant cohomology. In the fourth subsection we present some examples and we discuss the crucial property of reduction to the maximal torus of the cohomology, which now on twisted models plays an even more important role. Finally in the fifth subsection we reinterpret the proposed models as an example of a more general strategy which could be applied to a larger class of deformations.

\subsection{Classical models and Weil algebra}

We recall the classical construction of equivariant cohomology for the action of compact Lie group
$G$ on a smooth manifold $\mfdm$. The theory was originally formulated by Cartan \cite{car}; for a modern treatment a excellent reference is \cite{gs}.

One looks for a definition of equivariant cohomology $H_G(\mfdm)$ which is well defined for general actions, but that reduces to $H(\mfdm/G)$ for free actions. Since we expect $H_G(\mfdm)$ to satisfy homotopy invariance,
the idea is to deform $\mfdm$ into a homotopically equivalent space $\mfdm^{\prime}$ where the action is now free, and define $H_G(\mfdm)=H(\mfdm^{\prime}/G)$. A possible way is to consider a contractible space $E$ on which $G$ acts freely, so that we can put $\mfdm^{\prime}=\mfdm\times E$; of course at the end we have to prove that the definition does not depend on the choice of $E$. 

A natural choice for a $E$ is the total space of the universal $G$ bundle $ G\hookrightarrow EG \rightarrow BG $; we denote $X_G=(X\times EG)/G$. This leads to the following definition of equivariant cohomology, known as the Borel model.

\begin{df}
\label{bordef}
The equivariant cohomology of a smooth manifold $\mfdm$ acted upon by a compact Lie group $G$ is defined as the ordinary cohomology of the space $\mfdm_G$:
\begin{equation}
\label{borhg}
H_G(\mfdm) := H(\mfdm_G) = H((\mfdm\times EG)/G)
\end{equation}
where $EG$ is the total space of the universal $G$-bundle.
\end{df}

The problem with this definition is that $EG$ is finite dimensional only for $G$ discrete. A good recipe to overcome this problem is to find a finitely generated algebraic model for the algebra of differential forms over $EG$; this is where the Weil algebra comes into play.

\begin{df}
\label{koszdef}
The Koszul complex of a $n$-dimensional vector space $V$ is the tensor product between the symmetric and the exterior algebra of $V$
$$ \mathcal{K}_V = Sym (V) \otimes \wedge (V) $$
We assign to each element of $\bigwedge (V)$ its exterior degree, and to each element in $Sym^k(V)$ degree $2k$. The Koszul differential
$\kd$ is defined on generators
\begin{equation}
\label{kddef}
\kd (v\otimes 1) = 0 \qquad\qquad\qquad \kd (1\otimes v)= v\otimes 1
\end{equation}
and then extended as a derivation on the whole $\mathcal{K}(V)$.
\end{df}

A standard result (see e.g. \cite{gs} for a proof) is that the Koszul complex is acyclic, i.e. its cohomology reduces to the degree zero where it equals the scalar field.  

\begin{df}
\label{wadef}
The Weil algebra associated to a Lie group $G$ is the Koszul complex of $\g^{\ast}$, the dual of the Lie algebra of $G$. 
\end{df}

\begin{df}
\label{classkoszgen}
Let $\{e_a\}$ be a basis for $\g$. The set of Koszul generators of $W_{\g}$ is given by
\begin{equation}
\label{kosz}
e^a = e^a\otimes 1 \qquad\qquad\qquad \vartheta^a = 1\otimes e^a
\end{equation}
\end{df}

We are interested in the $\gt$-da structure of $W_{\g}$, i.e. the definition of operators $(L,i,d)$ on it. 

\begin{df}
\label{lionkdf}
The Lie derivative $L_a$ is defined by the coadjoint action of $\g$ on $\g^{\ast}$; on Koszul generators it reads
\begin{equation}
\label{lonk}
L_a(e^b) = - f_{ac}^{\phantom{ac}b} e^c \qquad\qquad\qquad L_a(\vartheta^b) = - f_{ac}^{\phantom{ac}b} \vartheta^c
\end{equation}
The interior derivative $i_a$ is given by
\begin{equation}
\label{ionk}
i_a(e^b) = - f_{ac}^{\phantom{ac}b} \vartheta^c \qquad\qquad i_a(\vartheta^b)=\delta_a^b
\end{equation}
The differential is the Koszul one; we then have $d_W(e^a)=0$ and $d_W(\vartheta^a)=e^a$.
\end{df}

These operators are extended by a (graded) Leibniz rule to the whole Weil algebra. Note that $L$ is of degree zero, $i$ of degree $-1$, $d_W$ of degree $1$ and the usual commutation relations among $(L,i,d)$ are satisfied.

A different set of generators for $W_{\g}$ is obtained by using horizontal (i.e. annihilated by interior derivatives) even elements.

\begin{df}
The set of horizontal generators for $W_{\g}$ is $\{u^a,\vartheta^a\}$ where
\begin{equation}
\label{horgen}
u^a := e^a + \frac{1}{2} f_{bc}^{\phantom{bc}a}\vartheta^b\vartheta^c
\end{equation}
\end{df}

With basic computations one can find the action of $(L,i,d)$ on horizontal generators; the new expressions are
\begin{equation}
\label{lidonu}
\begin{array}{clccl}
L_a(u^b) & = - f_{ac}^{\phantom{ac}b} u^c           & \qquad\qquad & i_a(u^b) & =0 \\
d_W(u^a) & = - f_{bc}^{\phantom{bc}a}\vartheta^bu^c & \qquad\qquad & d_W(\vartheta^a) & = u^a - \frac{1}{2}f_{bc}^{\phantom{bc}a}\vartheta^b\vartheta^c
\end{array}
\end{equation}
so that even generators are killed by interior derivative, hence the name horizontal.

Given a commutative $\gt$-da $A$ the tensor product $W_{\g}\otimes A$ is again a $\gt$-da with $L^{(tot)}=L\otimes 1 + 1\otimes L$ and the same rule for $i$ and $d$; this comes from the tensor structure of the category of $\env{\gt}$-module algebras. The basic subcomplex of a $\gt$-da is the intersection between invariant and horizontal elements. We have now all the ingredients to define the Weil model for equivariant cohomology.

\begin{df}
\label{weilmodeldef}
The Weil model for the equivariant cohomology of a commutative $\gt$-da $A$ is the cohomology of the basic subcomplex of
$W_{\g}\otimes A$:
\begin{equation}
\label{weilmodel}
H_G(A) = \left( (W_{\g}\otimes A)^G_{hor} , \; \delta = d_W\otimes 1 + 1\otimes d \right)
\end{equation}
\end{df}

The Weil model is the algebraic analogue of the Borel model with $A=\Omega(X)$,
$W_{\g}$ playing the role of differential forms on $EG$ and the basic subcomplex representing differential forms on the quotient space for free actions. A rigorous proof that topological and algebraic definitions are equivalent, a result known as the 'Equivariant de Rham Theorem', may be found for example in \cite{gs}. \\

Another well known algebraic model for equivariant cohomology of $\gt$-da's is the Cartan model; it defines equivariant cohomology as the cohomology of equivariant differential forms with respect to a 'completion' of the de Rham differential. We derive it as the image of an automorphism of the Weil complex $W_{\g}\otimes A$; the automorphism is usually referred as the Kalkman map \cite{kal} and is defined as

\begin{equation}
\label{kalkmap}
\phi = \mbox{exp }\{\vartheta^a\otimes i_a\}: W_{\g}\otimes A \longrightarrow W_{\g}\otimes A
\end{equation}

The image via $\phi$ of the basic subcomplex of $W_{\g}\otimes A$, the relevant part for equivariant cohomology, is easily described. 

\begin{prop}
The Kalkman map $\phi$ realizes an algebra isomorphism
\begin{equation}
\label{kalkiso}
(W_{\g}\otimes A)_{hor}^G\stackrel{\phi}{\simeq} \left(Sym(\g^{\ast})\otimes A\right)^G
\end{equation}
\end{prop}

The proof is obtained by direct computation; see \cite{kal} or \cite{gs}. The algebra $\left(Sym(\g^{\ast})\otimes A\right)^G$ appearing in $(\ref{kalkiso})$ will define
the Cartan complex and it is denoted by $\car{A}$. The differential on $\car{A}$ is induced from $\delta$ by the Kalkman map.

\begin{prop} The Cartan differential $d_G=\phi\, \delta_{|bas} \phi^{-1}$ on $\car{A}$ takes the form
\begin{equation}
\label{dg}
d_G=1\otimes d - u^a\otimes i_a
\end{equation}
\end{prop}

Again this can be proved by direct computation; we refer to \cite{kal}\cite{gs} for the details.

\begin{df}
\label{carmoddef}
The Cartan model for the equivariant cohomology of a commutative $\gt$-da $A$ is the cohomology
of the Cartan complex $\car{A}$:
\begin{equation}
\label{cartmod}
H_G(A)=\left( (Sym(\g^{\ast})\otimes A)^G, \, d_G=1\otimes d - u^a\otimes i_a \,\right)
\end{equation}
\end{df}

We make here a remark on the relation between Weil, Cartan and BRST differentials \cite{kal}. Denote by $M_W$ the differential algebra $W_G\otimes A$ with $\delta = d_W\otimes 1 + 1\otimes d$; it is possible to define another differential on the same algebra, the BRST operator
\begin{equation}
\label{brstdiff}
\delta_{BRST}= \delta + \vartheta^a\otimes L_a - u^a\otimes i_a
\end{equation}
We call $M_{BRST}$ the differential algebra $(W_{\g}\otimes A,\delta^{BRST})$; for the physical interpretation of $M_{BRST}$ see \cite{kal}. The Kalkman map is a $\gt$-da isomorphism from $M_W$ to $M_{BRST}$, i.e. it intertwines the two $\gt$-da structures. When restricted to $(W_M)_{|bas}$ its image is the Cartan model, now seen as the $G$-invariant subcomplex of the BRST model $M_{BRST}$; then also the Cartan differential $d_G$ is nothing but the restriction to the invariant subcomplex of the BRST differential $\delta_{BRST}$. We will show that it is possible to deform all the three models to the nc setting and keep the same relation among them; we wish to point out that this could be an interesting first step toward a definition of a nc BRST cohomology, with possible applications to nc gauge theories.

We end the section by noting that any homomorphism of $\gt$-da induces by functoriality a homomorphism between the corresponding equivariant cohomologies. For every $\gt$-da $A$ by taking the natural homomorphism $\mathbb{C}\rightarrow A$ we get a $H_G(\mathbb{C})=(Sym(\g^{\ast}))^G$ module structure on $H_G(A)$; the differential $d_G$ commutes with this module structure. $H_G(\mathbb{C})$ is called the basic cohomology ring.

\subsection{The noncommutative equivariant cohomology of Alekseev and Meinrenken}

In the previous section we introduced the Weil algebra as a finitely generated algebraic model for differential forms over $EG$. In the spirit of nc geometry an even more appropriate way to think of $W_{\g}$ is as the universal locally free object in the category of commutative $\gt$-da's $A$\cite{am2}. Indeed by using this approach we have a natural ways to define Weil algebras even in categories of deformed or nc $\gt$-differential algebras.

The first example of this strategy is the nc Weil algebra $\ncwa$ of Alekseev and Meinrenken \cite{am1}\cite{am2}, which they use to define equivariant cohomology in the category of nc $\gt$-da's. We will review their construction, and in the next sectiosn we will move to the category of twisted $\gt$-da's. A more detailed discussion on universal properties of these deformed Weil algebras is postponed to section $2.5$ and to a forthcoming paper \cite{lu2}. \\

The nc Weil algebra of \cite{am1} has a better formulation if we make an additional hypothesis: we demand that $\g$ is a quadratic Lie algebra, i.e. a Lie algebra carrying a nondegenerate $ad$-invariant quadratic form $B$ which can be used to canonically identify $\g$ with $\g^{\ast}$. The most natural examples of quadratic Lie algebras are given by semisimple Lie algebras, taking the Killing forms as $B$; since we already decided to restrict our attention to semisimple Lie algebras $\g$ in order to have more explicit expressions for the Drinfeld twists, this additional hypothesis fits well in our setting and we shall use it from now on.

\begin{df}
Let $(\g,B)$ be a quadratic Lie algebra. Fix a basis $\{e_a\}$ for $\g$ and let  $f_{ab}^{\phantom{ab}c}$ be the structure constants for this basis.
The super Lie algebra $\gtb$ is defined as the super vector space $\g^{(ev)}\oplus\g^{(odd)}\oplus\mathbb{C}\mathfrak{c}$,
with basis given by even elements $\{e_a, \mathfrak{c}\}$ and odd ones $\{\xi_a\}$, and brackets given  by
\begin{equation}
\label{brsg}
\begin{array}{clclcl}
\phantom{a}[e_a,e_b] & = f_{ab}^{\phantom{ab}c} \qquad  & [e_a, \xi_b] & = f_{ab}^{\phantom{ab}c}\xi_c \qquad & [\xi_a, \xi_b] & =B_{ab}\mathfrak{c} \\
\phantom{a}[e_a, \mathfrak{c}] & = 0 \qquad  & [\xi_a, \mathfrak{c}] & = 0 & &
\end{array}
\end{equation}
\end{df}

Using $\gtb$ the nc Weil algebra of \cite{am1} may be defined as (the quotient of) a super-enveloping algebra; this apparently trivial fact (not even explicitly stated in \cite{am1}) will be crucial in the following to realize a deformed $\gt$-da structure suitable for the nc setting.

\begin{df}
For quadratic Lie algebras $(\g,B)$ the noncommutative Weil algebra $\ncwa$ is defined as
\begin{equation}
\label{ncwaqq}
\ncwa=\env{\gtb}/ \langle\mathfrak{c}-1\rangle \simeq  \env{\g}\otimes Cl(\g,B).
\end{equation}
\end{df}

From now on we shall consider $\ncwa$ as a super enveloping algebra; formally we are working in $\env{\gtb}$ assuming implicitly every time $\mathfrak{c}=1$. Moreover the decomposition of $\ncwa$ in the even part $\env{\g}$ and an odd part $Cl(\g,B)$ is by the time being only true as a vector space isomorphism; to become an algebra isomorphism we have to pass to even generators which commute
with odd ones: this will be done below.

We are interested in the $\gt$-da structure of $\ncwa$. The main difference with the classical Weil algebra
is that the action of $(L,i,d)$ may now be realized by inner derivations.

\begin{df}
\label{gdancwa}
On a generic element $X\in\ncwa$ the actions of $L$ and $i$ are given by
\begin{equation}
\label{il}
L_a(X):=ad_{e_a}(X) \qquad \qquad i_a(X):= ad_{\xi_a}(X)
\end{equation}
On generators one has
\begin{equation}
\label{ilongen}
\begin{array}{lcl}
L_a(e_b)=[e_a,e_b]= f_{ab}^{\phantom{ab}c} e_c       & \qquad \qquad \qquad & i_a(e_b)= [\xi_a,e_b]= f_{ab}^{\phantom{ab}c}\xi_c \\
L_a(\xi_b)=[e_a,\xi_b]= f_{ab}^{\phantom{ab}c} \xi_c & \qquad \qquad \qquad & i_a(\xi_b)= [\xi_a,\xi_b] = B_{ab}\mathfrak{c}
\end{array}
\end{equation}
\end{df}

Thus $L_a$ and $i_a$ are derivations (thanks to the primitive coproduct of $e_a$ and $\xi_a$ in $\env{\gt}$) and their action agrees with the commutator of $e_a$ and $\xi_a$ in $\ncwa$.

\begin{df}
The differential $\nwd$ on the noncommutative Weil algebra $\ncwa$ is the Koszul differential
$\nwd(e_a)=0$ , $\nwd(\xi_a)=e_a$, so that $(\ncwa , \nwd)$ is an acyclic differential algebra.
\end{df}

Following the terminology Def \ref{koszdef} the set of generators $\{e_a,\xi_a\}$ of $\ncwa$ will be called of Koszul type. It is often more convenient to use horizontal generators. These are introduced by the transformation
\begin{equation}
\label{u}
u_a:= e_a+ \frac{1}{2}f_a^{\phantom{a}bc}\xi_b\xi_c
\end{equation}
where we use $B$ to raise and lower indices. One can easily verify that $\{u_a, \xi_a\}$ is another set of generators for $\ncwa$, with relations (compare with (\ref{brsg})):
\begin{equation}
\label{relucl2}
[u_a,u_b]=f_{ab}^{\phantom{ab}c} u_c \qquad\qquad [u_a,\xi_b]=0   \qquad\qquad [\xi_a, \xi_b]=B_{ab}
\end{equation}
Note that $u_a$ generators realize the same Lie algebra $\g$ of $\{e_a\}$, but now decoupled from the odd part, so that using these generators we can write $\ncwa \simeq \env{\g}\otimes Cl(\g,B)$ as an algebra isomorphism. We skip the proof of the following elementary restatement of relations in Def \ref{gdancwa}.

\begin{prop}
The $\gt$-da structure, still given by adjoint action of generators $\{e_a, \xi_a\}$, now on $\{u_a, \xi_a\}$ reads:
\begin{equation}
\label{ilongen2}
\begin{array}{rclclcl}
L_a(u_b) & = & f_{ab}^{\phantom{ab}c} u_c                 & \qquad \qquad \qquad & L_a(\xi_b) & = & f_{ab}^c \xi_c \\
i_a(u_b) & = & 0                                          & \qquad \qquad \qquad & i_a(\xi_b) & = & B_{ab}\\
\nwd(u_a)& = & -f_a^{\phantom{a}bc}\xi_b u_c              & \qquad \qquad \qquad & \nwd(\xi_a)& = & u_a - \frac{1}{2}f_a^{\phantom{a}bc}\xi_b\xi_c
\end{array}
\end{equation}
\end{prop}

The operator $\nwd$ may be expressed as an inner derivation as well: indeed it is given by the commutator with an element $\mathcal{D}\in (\ncwa^{(3)})^G$. There are several ways (depending on the choice of generators used) one can write $\mathcal{D}$, and the simplest one for our calculations is
\begin{equation}
\label{D}
\mathcal{D}=\frac{1}{3}\,\xi^ae_a+\frac{2}{3}\,\xi^au_a
\end{equation}
For a generic element $X\in\ncwa$ we can then write $\nwd(X)=[\mathcal{D},X]$. Notice that $\ncwa$ is a filtered differential algebra, with associated graded differential algebra the classical Weil algebra $W_{\g}$; the $\gt$-da structure of $\ncwa$ agrees with the classical one if we pass to $Gr(\ncwa)$. 

Given any $\gt$-da $A$ the tensor product $\ncwa\otimes A$ gets a natural $\gt$-da structure (which is unbraided since we are still considering $\env{\gt}$-module algebras). Following the classical construction we define equivariant cohomology as the cohomology of the basic subcomplex of $W_{\g}\otimes A$.

\begin{df}\cite{am1}
\label{ncwm}
The Weil model for the equivariant cohomology of a nc $\gt$-differential algebra $A$ is the cohomology of the complex
\begin{equation}
\mathcal{H}_G(A) = \left( (\ncwa\otimes A)^G_{(hor)} , \, \delta^{(tot)} = \nwd\otimes 1 + 1\otimes d \right)
\end{equation}
\end{df}

There nc analogue of the Kalkman map (\ref{kalkmap}), expressed using generators of $\ncwa$, is
\begin{equation}
\label{nckalkmap}
\Phi = \mbox{exp } \{ \xi^a\otimes i_a\} : \ncwa\otimes A \longrightarrow \ncwa\otimes A
\end{equation}
By a proof completely similar to the classical one, in \cite{am1} it is shown how $\Phi$ intertwines the action of $L^{(tot)}$ and $i^{(tot)}$, leading to the following result.

\begin{prop}
The nc Kalkman map $\Phi$ defines a vector space isomorphism
\begin{equation}
(\ncwa\otimes A)_{hor}^G \stackrel{\Phi}{\simeq} (\env{\g}\otimes A)^G
\end{equation}
\end{prop}

The main difference between the classical and the nc Kalkman map is that $\xi^a\otimes i_a$ is no longer a derivation; for this reason $\Phi$ is not an algebra homomorphism, and the natural algebra structure on $(\env{\g}\otimes A)^G$ does not agree with the one induced by $\Phi$. Before looking at the algebra structure of the image of the Kalkman map we describe the induced differential.

\begin{prop}\cite{am1}
\label{ncdcartthm}
The nc Cartan differential $d_G$ induced from $\delta^{(tot)}=\dw\otimes 1 + 1\otimes d$ by the Kalkman map $\Phi$ via
$d_G=\Phi(\nwd\otimes 1 + 1\otimes d)_{|bas}\Phi^{-1}$ takes the following expression
\begin{equation}
\label{ncdcart}
d_G = 1\otimes d - \frac{1}{2} (u^a_{(L)} + u^a_{(R)}) \otimes i_a + \frac{1}{24} f^{abc} (1\otimes i_ai_bi_c)
\end{equation}
where with $u^a_{(L)}$ (resp. $u^a_{(R)}$) we denote the left (resp. right) multiplication for $u^a$.
In particular $d_G$ commutes with $L$ and squares to zero on $(\env{\g}\otimes A)^G$.
\end{prop}

As previously discussed, the Kalkman map is a $\gt$-da iso between the Weil model and the BRST model; when we restrict the image of $\Phi$ to the basic subcomplex we find the Cartan model \cite{kal}. We can then interpret the image of the nc Kalkman map as a nc BRST model; with a direct computation one can check that the nc BRST differential is 
\begin{equation}
\label{ncbrstd}
\delta_{BRST}=\Phi(\nwd\otimes 1+1\otimes d)\Phi^{-1}= d_G+ \nwd\otimes 1 + \xi^a\otimes L_a
\end{equation}
where by $d_G$ we mean the nc Cartan differential (\ref{ncdcart}); note that as expected $(\delta_{BRST})_{|bas}=d_G$. We denote the complex $(\nccarl{A}, d_G)$ by $\nccar{A}$. Its ring structure is induced by the Kalkman map; by definition on $u_i\otimes a_i \in \nccarl{A}$ we have
\begin{equation}
\label{moltcart}
(u_1\otimes a_1)\odot (u_2\otimes a_2) := \Phi \left ( \Phi^{-1} (u_1\otimes a_1)\cdot_{\ncwa\otimes A}
\Phi^{-1} (u_2\otimes a_2) \right)
\end{equation}

\begin{prop}\cite{am1}
\label{ringcartan}
The ring structure of $\nccar{A}$ defined in (\ref{moltcart}) takes the explicit form
\begin{equation}
\label{nccm}
(u_1\otimes a_1)\odot(u_2\otimes a_2) = (u_1u_2)\otimes \cdot_{A} \left(\mbox{exp }\{B^{rs}i_r\otimes i_s\}(a_1\otimes a_2)\right)
\end{equation}
Note that $d_G$ is a derivation of $\odot$.
\end{prop}

\begin{df}
\label{nccartanmodeldef}
The Cartan model for the equivariant cohomology of a nc $\gt$-da $A$ is the cohomology of the complex $(\nccar{A},d_G)$:
\begin{equation}
\mathcal{H}_G(A) = \left( \nccarl{A} , \,
d_G = 1\otimes d - \frac{1}{2} (u^a_{(L)} + u^a_{(R)}) \otimes i_a + \frac{1}{24} f^{abc}\otimes i_ai_bi_c \right)
\end{equation}
The ring structure $\odot$ of $\nccar{A}$ is given in (\ref{nccm}).
\end{df}

Note that for abelian groups the Cartan model reduces to the classical one; in the non abelian case this ring structure is not compatible with a possible pre-existing grading on $A$. The only structure left in $\nccar{A}$ is a double filtration; its associated graded differential module is a double graded differential model and agrees with the classical Cartan model.

We finally stress that these nc Weil and Cartan model do apply to nc algebras, but the request is that 
the $\gt$-da structure is undeformed. We are rather interested in nc algebras where the noncommutatitivy is strictly related to a deformed $\gt$-da structure; basically we are interested in a different category, so we need different models. 

\subsection{Twisted noncommutative equivariant cohomology}

In this section we introduce models for the equivariant cohomology of twisted nc $\gt$-da's $A_{\chi}$, i.e. $\gt$-da's deformed by a Drinfel'd twist as in Thm \ref{covdef}. We show how it is possible to mimic the construction of Alekseev and Meinrenken of the previous subsection. Basically we deform the nc Weil algebra $\ncwa$ using the same $\chi$ which realizes the deformation of the nc algebra $A_{\chi}$; we keep considering quadratic Lie algebras $\g$, so that $\ncwa$ is an enveloping algebra and the twsist $\chi$ acts naturally on it. The definition of Weil and Cartan models will follow as usual from the cohomology of the appropriate subcomplexes.

The construction we are going to present works for arbitrary twisted $\gt$-da's, even in the cases where the form of the twist element $\chi$ is unknown. Obviously if one wants to deal with explicit expressions and computations, like the ones presented here, an explicit form of $\chi$ is crucial; in what follows we will continue to use the Drinfel'd twist $\chi$ given in (\ref{twistel}).

The Weil algebra of the category of twisted nc $\gt$-da's will have a twisted $\gt$-da structure; a natural candidate is the Drinfel'd twist of the nc Weil algebra $\ncwa$.

\begin{df}
\label{tncwadef}
Let $\g$ be a quadratic Lie algebra, and $A_{\chi}$ a twisted nc $\gt$-da. The twisted nc Weil algebra $\tncwa$ is defined as the Drinfeld twist of $\ncwa$ by the same $\chi$, now viewed as
an element in $\left(\ncwa\otimes\ncwa\right)^{(ev)}_{[[\theta]]}$ .
\end{df}

The generators of the twist $\chi$ need not to belong to $\g$; we already discussed the same fact when deforming actions on nc algebras in section $1.3$; we recall that in that case the relevant Lie algebra is the product between $\ttt$, the torus which contains the generators of $\chi$, and $\g$, the symmetry whose action is relevant for equivariant cohomology. Of course the interesting case is when $\ttt$ and $\g$ do not commute, otherwise the twist is trivial. In what follows we will directly assume that $\g$ contains the generators of the twist.

We want to describe the $\gt$-da structure of the twisted Weil algebra. Following the usual notation we denote even and odd generators of $\tncwa$ by $\{e_i, e_r, \xi_i, \xi_r\}$, distinguishing between Cartan (index $i$) and root (index $r$) elements of $\gt$. We already computed the twisted coproduct of the even subalgebra (see Prop \ref{twcopr}) and of odd generators $\xi_a$ (see Prop \ref{twodd}).
Recall also that, as showed in Prop \ref{twantip}, for this class of Drinfeld twist elements $\chi$ the antipode is undeformed.

The $\gt$-da structure of the nc Weil algebra $\ncwa$ has been realized by the adjoint action with respect to even generators (for the Lie derivative), odd generators (for the interior derivative) and by commutation with a fixed element in the center (for the differential). We use the same approach for $\tncwa$, the only difference is that now from the general formula for the adjoint action on a super Hopf algebra
\begin{equation}
ad_{Y}(X)= \sum (-1)^{|X||(Y)_{(2)}|} (Y)_{(1)}X (S(Y)_{(2)})
\end{equation}
we see that the twisted coproduct generates a twisted adjoint action even on single generators.

\begin{df}
The action of $L$ and $i$ on $\tncwa = \envt{\gt}$ is given by the adjoint action with respect to even and odd generators. In particular $L_i=ad_{e_i}$ and $i_i=ad_{\xi_i}$ are the same as in the untwisted case. On the contrary for roots elements the operators $L_r$ and $i_r$ are modified even on single generators:
\begin{equation}
\label{litwe}
\begin{split}
L_r(X) & = ad^{\chi}_{e_r}(X) = e_rX\lambda_r - \lambda_r X e_r\\
i_r(X) & = ad^{\chi}_{\xi_r}(X) = \xi_r X \lambda_r + (-1)^{|X|} \lambda_r X \xi_r
\end{split}
\end{equation}
\end{df}
Expressing explicitly this action on $\{e_a, \xi_a\}$ we have (one should compare with (\ref{ilongen})):
\begin{equation}
\label{longentw}
\begin{array}{rlccrl}
L_j(e_a)= & f_{ja}^{\phantom{ja}b}\,e_b        & & &   L_j(\xi_a)= & f_{ja}^{\phantom{ja}b}\,\xi_b \\
L_r(e_i)= & e_re_i\lambda_r - \lambda_re_ie_r  & & &   L_r(\xi_i)= & e_r\xi_i\lambda_r-\lambda_r\xi_ie_r \\
        = & -r_i\lambda_re_r                   & & &             = & -r_i\lambda_r\xi_r \\
L_r(e_s)= & e_re_s\lambda_r - \lambda_re_se_r  & & &   L_r(\xi_s)= & e_r\xi_s\lambda_r - \lambda_r\xi_se_r
\end{array}
\end{equation}
\begin{equation}
\label{iongentw}
\begin{array}{rlccrl}
i_j(e_a)= & f_{ja}^{\phantom{ja}b}\,\xi_b         & & & i_j(\xi_a) = & B_{ja}=\delta_{ja}\\
i_r(e_i)= & \xi_re_i\lambda_r - \lambda_re_i\xi_r & & & i_r(\xi_i) = & \xi_r\xi_i\lambda_r+\lambda_r\xi_i\xi_r \\
        = & -r_i\lambda_r\xi_r                    & & &            = & \lambda_rB_{ri}=0 \\
i_r(e_s)= & \xi_re_s\lambda_r - \lambda_re_s\xi_r & & & i_r(\xi_s) = & \xi_r\xi_s\lambda_r + \lambda_r\xi_s\xi_r
\end{array}
\end{equation}
where we use $i, j$ for Cartan indexes, $r, s$ for roots indexes and $a,b$ for generic indexes.
On products one just applies the usual rule for the adjoint action
\begin{equation}
ad_{Y}(X_1X_2)= (ad_{Y_{(1)}}X_1)(ad_{Y_{(2)}}X_2)
\end{equation}
which shows that $L_r$ and $i_r$ are twisted derivations.

Due to the presence of the $\lambda_r$ terms the classical generators $\{e_a, \xi_a\}$ are no longer closed under the action of $L,i$. There is however another set of generators (we will call them quantum generators for their relation to quantum Lie algebras, see below) which is more natural.
\begin{df}
\label{qgdef}
The quantum generators of $\tncwa$ are 
\begin{equation}
\label{qg}
X_a:=\lambda_a e_a \qquad \qquad \eta_a:= \lambda_a\xi_a
\end{equation}
Recall from (\ref{lambda}) that for $a=i$ we have $\lambda_i=1$, so $X_i=e_i$. We define also coefficients
\begin{equation}
\label{q}
q_{rs}:=\mbox{exp }\{\frac{i}{2}\theta^{kl}r_ks_l\}
\end{equation}
with properties $q_{sr}=q_{rs}^{-1}$ and $q_{rs}=1$ if $r=-s$; we also set $q_{ab}=1$ if at least one index is of Cartan type (due to the vanishing of the correspondent root vector).
\end{df}
The following relations, easily proved by direct computation,  will be very useful:
\begin{equation}
\label{rel1}
\begin{array}{rlcrl}
\lambda_r\lambda_s & = \lambda_{r+s}      & \qquad & \lambda_r\lambda_s & = \lambda_s\lambda_r   \\
\lambda_re_s       & = q_{rs}e_s\lambda_r & \qquad & \lambda_r\xi_s     & = q_{rs}\xi_s\lambda_r \\
L_{\lambda_r}e_s   & = q_{rs}e_s          & \qquad & L_{\lambda_r}\xi_s & = q_{rs}\xi_r
\end{array}
\end{equation}
and since all $\lambda_r$'s commute with each other, the same equalities hold for $X_r$ and $\eta_r$.
Using the definition of the adjoint action, the previous relations (\ref{rel1}) and the commutation rules between $\{e_a, \xi_a\}$ in $\tncwa$ we can express by straightforward computations the twisted $\gt$-da structure on quantum generators.

\begin{prop}
The action of $L$ and $i$ on quantum generators $\{X_a,\eta_a\}$ of $\tncwa$ is
\begin{equation}
\label{lionXeta}
\begin{array}{rlccrl}
L_aX_b    & = f_{ab}^{\phantom{ab}c}X_c    & & & i_aX_b    & = f_{ab}^{\phantom{ab}c}\eta_c \\
L_a\eta_b & = f_{ab}^{\phantom{ab}c}\eta_c & & & i_a\eta_b & = B_{ab}
\end{array}
\end{equation}
\end{prop}

Note that this is exactly the same action we have in the classical case (\ref{ilongen}).
The difference however is that we keep acting on quantum generators with classical generators: $L_aX_b=ad_{e_a}X_b\neq ad_{X_a}X_b$. 

We make a quick digression on the meaning of quantum generators and their link with quantum Lie algebras, even if this is not directly related to the construction of equivariant cohomology.
The fact that the generators $\{e_a, \xi_a\}$ are not closed under the deformed adjoint action is a typical feature of quantum enveloping algebras $\mathfrak{U}_{q}(\g)$ where the deformation involves the Lie algebra structure of $\g$ (contrary to what Drinfeld twists do). Since $\g$ can be viewed as the closed $ad$-submodule of $\env{\g}$ one can try to recover a quantum Lie algebra inside $\mathfrak{U}_{q}(\g)$ by defining $\g_{q}$ as a closed $ad$-submodule of $\mathfrak{U}_{q}(\g)$ with quantum Lie bracket given by the adjoint action. The quantum Lie brackets are linear, $q$-skewsymmetric and satisfy a deformed Jacobi identity \cite{qla}.

In the Drinfel'd twist case the deformation of the coproduct in $\envt{\g}$ leads to a deformation of the adjoint action, but the brackets $[e_a,e_b]$ are unchanged; thus $ad_{e_r}(e_s)$ is no more equal to $[e_r,e_s]$. However $\{X_a\}$ are generators of a closed $ad$-submodule (see (\ref{rel1})), so we can define quantum Lie brackets $[\; , \, ]_{(\chi)}$ using the twisted adjoint action, obtaining a quantum Lie algebra structure $\g_{\chi}$:
\begin{equation}
\label{qlabr}
\begin{array}{ll}
\,[X_i,X_j]_{(\chi)}     & := ad^{\chi}_{X_i}X_j     = 0 \\
\,[X_i, X_r]_{(\chi)}    & := ad^{\chi}_{X_i}X_r     = r_iX_r = - [X_r, X_i]_{(\chi)}  \\
\,[X_{-r},X_r]_{(\chi)}  & := ad^{\chi}_{X_{-r}}X_r  = \sum r_iX_i = [X_r,X_{-r}]_{(\chi)} \\
\,[X_r,X_s]_{(\chi)}     & := ad^{\chi}_{X_r}X_s     = q_{rs}f_{rs}^{\phantom{rs}r+s}X_{r+s} \\
\,[X_s,X_r]_{(\chi)}     & := ad^{\chi}_{X_s}X_r     = q_{sr}f_{sr}^{\phantom{rs}r+s}X_{r+s} = -(q_{rs})^{-1}f_{rs}^{\phantom{rs}r+s}X_{r+s}
\end{array}
\end{equation}
The $q$-antisymmetry is explicit only in the $[X_r,X_s]_{(\chi)}$ brackets since $q_{ab}\neq 1$ if and only if both indexes are root type. The same result holds also for the odd part of $\bar{\gtb}$, so we may consider $\{X_a,\eta_a , \mathfrak{c}\}$ as a base for the quantum (super) Lie algebra inside $\envt{\bar{\gtb}}$. The last observation is that $\triangle_{\chi}X_r = X_r\otimes 1 + \lambda_r^2\otimes X_r$, so if we want $\g_{\chi}$ to be closed also under the coproduct, we may consider mixed generators $\{\Lambda_j,X_r\}$ where the Cartan-type generators are defined as group-like elements $\Lambda_j:= \mbox{exp }\{\frac{i}{2}\theta^{jl}H_l\}$. Now $\{\Lambda_j,X_r,\mathfrak{c}\}$ describe a different quantum Lie algebra $\g_{\chi}^{\prime}$, due to the presence of group-like elements; the structure of $\g_{\chi}$ is recovered taking the first order terms in $\theta$ of the commutators involving $\Lambda_j$'s.

We come back to equivariant cohomology and the twisted Weil algebra; it is useful to introduce horizontal generators. 

\begin{df}
\label{qhorg}
The quantum horizontal generators of $\tncwa$ are defined by
\begin{equation}
\label{k}
K_a:= \lambda_au_a= \lambda_a(e_a+\frac{1}{2}f_a^{\phantom{a}bc}\xi_b\xi_c)=X_a-\frac{1}{2}\eta^bad_{X_b}(\eta_a)
\end{equation}
\end{df}
They are indeed in the kernel of the twisted interior derivative
\begin{equation}
\label{iK}
  i_aK_b= ad^{\chi}_{\xi_a}(\lambda_bu_b)=\xi_a\lambda_bu_b\lambda_a-\lambda_a\lambda_bu_b\xi_a =0
\end{equation}
and their transformation under $L_a$ is given by
\begin{equation}
\label{lK}
L_aK_b= ad^{\chi}_{e_a}(\lambda_bu_b)=e_a\lambda_bu_b\lambda_a - \lambda_a\lambda_bu_be_a = f_{ab}^{\phantom{ab}c}K_c
\end{equation}

The last thing to describe is the action of the differential $\nwd$. Recall that in $\ncwa$ we had $\nwd(X)=[\mathcal{D},X]$, and this is still true in $\tncwa$. In fact
$$\mathcal{D}=\frac{1}{3}\,\xi^ae_a+\frac{2}{3}\,\xi^au_a=\frac{1}{3}\,\eta^aX_a+\frac{2}{3}\,\eta^aK_a$$
Moreover $\nwd$ being a commutator, the Jacobi identity assures it is an untwisted derivation. This is not surprising: the twisted $\gt$-da structure of an algebra does not change the action of the differential. Note that $\eta^a=\lambda_a^{-1}\xi^a$ and $\nwd\lambda_a=[\mathcal{D},\lambda_a]=0$. For even generators we have
\begin{equation}
\label{dK}
\dw(K_a)=\lambda_ad_W(u_a)=-f_a^{\phantom{a}bc}\lambda_a\xi_bu_c=
-f_a^{\phantom{a}bc}\lambda_b\lambda_c\xi_bu_c=-q_{ab}f_a^{\phantom{a}bc}\eta_bK_c
\end{equation}
where if we raise the index of $\eta$ we take in account the $\lambda$ inside $\eta$
\begin{equation}
-q_{ab}f_a^{\phantom{a}bc}\eta_bK_c = -q_{ba}f_{ab}^{\phantom{ab}c}\eta^bK_c
\end{equation}
For odd generators
\begin{equation}
\label{deta}
  \dw(\eta_a) =\lambda_ae_a= \lambda_a(u_a - \frac{1}{2} f_a^{\phantom{a}bc}\xi_b\xi_c) =
  K_a -\frac{1}{2} q_{ba}f_{ab}^{\phantom{ab}c}\eta^b\eta_b
\end{equation}

We have found all the relations which define a twisted $\gt$-da structure on $\tncwa$.
At this point we can define a Weil complex for any twisted $\gt$-da $A_{\chi}$; nc differential forms $\Omega(\mfdm_{\theta})$ provide a natural example to which the theory applies. The Weil complex involves the tensor product between the two twisted $\gt$-da's $\tncwa$ and $A_{\chi}$. We already showed that this construction depends on the quasitriangular structure of $\envt{\gt}$ (see Prop \ref{braidtens}); the deformed $\mathcal{R}$ matrix is (see Thm \ref{twistthm})
$\mathcal{R}^{\chi} = \chi_{21}\mathcal{R}\chi^{(-1)}$ (with $\chi_{21}= \chi^{(2)}\otimes\chi^{(1)}$).
Since the original $\mathcal{R}$ matrix of $\mathfrak{U}(\gt)$ is trivial we have the simple expression
\begin{equation}
\mathcal{R}^{\chi}=\chi^{-2} = \mbox{ exp}\{ i\theta^{kl}H_k\otimes H_l\}
\end{equation}
We introduce the twisted nc Weil model; the relevant difference is that now the tensor product between $\tncwa$ and $A_{\chi}$ is in the braided monoidal category of $\envt{\gt}$-module algebras.

\begin{df}
\label{tncwm}
The Weil model for the equivariant cohomology of a twisted nc $\gt$-da $A_{\chi}$ is the cohomology of the complex
\begin{equation}
  \mathcal{H}^{\chi}_G(A_{\chi})=\left( (\tncwa\widehat{\otimes}A_{\chi})_{bas}, \, \delta = \dw\otimes 1 + 1\otimes d \,\right)
\end{equation}
\end{df}
The basic subcomplex is taken with respect to $L^{tot}$ and $i^{tot}$; these operators act on $\tncwa\widehat{\otimes}A_{\chi}$ with the covariant rule $L^{tot}_X= L_{X_{(1)}}\otimes L_{X_{(2)}}$ using the twisted coproduct. We can use the $G$-invariance to explicitly compute the effect of the braiding on the multiplicative structure of the Weil model.
\begin{prop}
\label{multtncwm}
Let $A_{\chi}$ be a twisted nc $\gt$-da, with $A$ a graded-commutative $\gt$-da.
The multiplication in the Weil complex $(\tncwa\widehat{\otimes}A_{\chi})_{bas}$, according to the general formula (\ref{braidmult}), reads
\begin{equation}
\label{brpr}
(u_1\otimes\nu_1)\cdot (u_2\otimes b_2) = (-1)^{|\nu_1||\nu_2|} u_1u_2\otimes \nu_2\cdot_{\chi}\nu_1
\end{equation}
\end{prop}

\noindent\underline{\textsf{Proof}}: By direct computation, applying Lemma \ref{l1} to the left hand side and using $G$-invariance:
\begin{equation*}
\begin{split}
\sum_n & \frac{(i\theta^{\alpha\beta})^n}{n!}  \, u_1(H_{\beta}^n u_2) \otimes (H_{\alpha}^n \nu_1)\cdot_{\chi} \nu_2 = \\
 & = \sum_n \frac{(-i\theta^{\alpha\beta})^n}{n!} u_1u_2\otimes (H_{\alpha}^n\nu_1)\cdot_{\chi} (H_{\beta}^n\nu_2) =
 u_1u_2 \otimes \cdot (\chi^2\chi^{-1}\triangleright \nu_1\otimes\nu_2) = \\
 & = u_1u_2\otimes \cdot(\chi\triangleright \nu_1\otimes\nu_2) = (-1)^{|\nu_1||\nu_2|} u_1u_2\otimes \nu_2\cdot_{\chi}\nu_1
\qquad\qquad _{\blacksquare}
\end{split}
\end{equation*}

We want to compare $(\tncwa\widehat{\otimes}A_{\chi})_{bas}$ with the Weil complex of \cite{am1}. According to the philosophy of Drinfel'd twist deformations, namely to preserve the vector space structure and to deform only the algebra structure of $\gt$-da's, we find that they are isomorphic roughly speaking as 'vector spaces'; the precise statement, since we are comparing quantities depending on formal series in $\theta$, involves topologically free $\mathbb{C}_{[[\theta]]}$ modules, or $\theta$-adic vector spaces.

\begin{prop}
\label{isobastfm}
There is an isomorphism of (graded) topologically free $\mathbb{C}_{[[\theta]]}$ modules
$$(\tncwa\widehat{\otimes}A_{\chi})_{bas}\simeq \left((\ncwa\otimes A)_{bas}\right)_{[[\theta]]}$$
\end{prop}

\proof{ We first show the inclusion $\cwm\subseteq\twm$. Take
$$u\otimes\nu\in\cwm \Rightarrow (L\otimes 1+1\otimes L)(u\otimes\nu)=0$$
The $\g$ invariance property applied to powers of toric generators gives
$$ H_{\alpha}^n u \otimes \nu = (-1)^n u \otimes H_{\alpha}^n \nu$$
and in particular $\lambda_r u\otimes \nu = u\otimes \lambda_r^{-1} \nu $. This can be used to compute
\begin{equation*}
(L_r\otimes\lambda_r^{-1} + \lambda_r\otimes L_r)(u\otimes\nu) = (L_r\lambda_r\otimes 1 - \lambda_rL_r\otimes 1)(u\otimes\nu)
 = ([L_r,\lambda_r]\otimes 1)(u\otimes\nu) = 0
\end{equation*}
A similar short calculation (just writing $i_r$ instead of $L_r$) gives the analogous result for $i_r$ as well; so we showed that $u\otimes\nu\in\twm$. For the opposite inclusion, take now $v\otimes\eta \in \twm$; this implies 
$$(L_r\otimes\lambda^{-1}_r+\lambda_r\otimes L_r)(v\otimes\eta)=0$$ 
and in particular again $\lambda_r v\otimes \eta = v\otimes \lambda_r^{-1}\eta$. We use these two equalities to compute
\begin{equation*}
L_r v \otimes \eta = L_r\lambda_r^{-1} v\otimes \lambda^{-1}_r \eta = - (1\otimes L_r\lambda_r)(1\otimes\lambda_r)(v\otimes\eta) = - v\otimes L_r\eta
\end{equation*}
Substituting again $L_r$ with $i_r$ we easily find the same result for $i_r$, and this proves that $v\otimes\eta \in \cwm$. The linearity of the operators with respect to formal series in $\theta$ and the compatibility of the eventual grading (coming from $A$) with the $\mathbb{C}_{[[\theta]]}$-module structure complete the proof.}

The previous result easily generalizes to the associated equivariant cohomologies, since the differentials for both the complexes are the same.

\begin{prop}
\label{propisoam}
There is an isomorphism of (graded) topologically free modules
\begin{equation}
\label{isocohomtfm}
\mathcal{H}^{\chi}_G(A_{\chi}) \simeq \mathcal{H}_G(A)_{[[\theta]]}
\end{equation}
\end{prop}

\proof{ Since both $\mathcal{H}^{\chi}_G(A_{\chi})$ and  $\left( \mathcal{H}_G(A) \right)_{[[\theta]]}$ are defined starting from the respective basic subcomplexes with the same $\mathbb{C}_{[[\theta]]}$-linear differential $\delta = \nwd \otimes + 1\otimes d$ the isomorphism of Prop \ref{isobastfm} lifts to the cohomologies.}

Roughly speaking we are saying that our twisted equivariant cohomology is equal to the trivial formal series extension of the nc cohomology of Alekseev and Meinrenken, as 'vector space' over $\mathbb{C}_{[[\theta]]}$ (i.e. as topologically free $\mathbb{C}_{[[\theta]]}$-module). This is not surprising, since we expect the deformation coming from the Drinfel'd twist to be visible only at the ring structure level.

We now pass to the construction of a twisted nc Cartan model. Basically we need to twist the nc Kalkman map of \cite{am1} in order to intertwine the twisted Lie and interior derivative which define the basic subcomplex.

\begin{df}
The twisted nc Kalkman map
$$\Phi^{\chi}:\tncwa\widehat{\otimes}A_{\chi}\rightarrow \tncwa\widehat{\otimes}A_{\chi}$$
is the conjugation by the twist element $\chi$ of the nc Kalkman map $\Phi$
\begin{equation}
\Phi^{\chi}= \chi\Phi\chi^{-1} \qquad\qquad \mbox{with } \, \Phi = \mbox{exp }\{\xi^a\otimes i_a\}
\end{equation}
\end{df}

\begin{prop}
There is an isomorphism of topological free $\mathbb{C}_{[[\theta]]}$-modules
\begin{equation}
\label{3w}
(\tncwa\widehat{\otimes}A_{\chi})_{bas} \stackrel{\Phi^{\chi}}{\simeq}
(\tncwa\widehat{\otimes}A_{\chi})^G_{i_a\otimes\lambda_a^{-1}} =
(\tncwa\widehat{\otimes}A_{\chi})^G_{i_a\otimes 1}
\end{equation}
\end{prop}

\proof{ First note that $\Phi^{\chi}$ is invertible with $(\Phi^{\chi})^{-1}=\chi\Phi^{-1}\chi^{-1}$. To prove equivariance of $\Phi^{\chi}$ note that the $\chi$ coming from the twisted coproduct cancels with the $\chi$ in $\Phi^{\chi}$:
\begin{equation*}
\begin{split}
\Phi^{\chi}L^{(tot)}_r(\Phi^{\chi})^{-1} & = (\chi\Phi\chi^{-1})(\chi\triangle (u_r) \chi^{-1})(\chi\Phi^{-1}\chi^{-1}) =
\chi (\Phi \triangle (u_r) \Phi^{-1}) \chi^{-1} = \\
 & = \chi \triangle (u_r) \chi^{-1} = L_r^{(tot)}
\end{split}
\end{equation*}
where we used the equivariance of $\Phi$ with respect to the untwisted $L^{(tot)}$. A similar computation for $i^{(tot)}$ gives
\begin{equation*}
\begin{split}
\Phi^{\chi}i^{(tot)}_r(\Phi^{\chi}){-1} & = (\chi\Phi\chi^{-1})(\chi\triangle (\xi_r) \chi^{-1})(\chi\Phi^{-1}\chi^{-1}) =
\chi (\Phi \triangle (\xi_r) \Phi^{-1}) \chi^{-1} = \\
 & = \chi (i_r\otimes 1) \chi^{-1} = i_r\otimes\lambda_r^{-1}
\end{split}
\end{equation*}
The last equality comes easily from the computation of $\chi(i_r\otimes 1)\chi^{-1}$ expanding $\chi$ at various orders in $\theta$. Finally we get the right hand side of $(\ref{3w})$ using $\lambda_a\otimes\lambda_a= 1\otimes 1$ on basic elements.}

In the untwisted setting we have $(\ncwa)_{hor}\simeq\env{\g}$. Here $(\tncwa)_{hor}=\{K_a\}\neq\envt{\g}$, that is the horizontal subalgebra of $\tncwa$ is spanned by quantum horizontal generators $K_a$ (see Def \ref{qhorg}) which do not describe any enveloping algebra. We will use the following notation to refer to the image of $\Phi^{\chi}$:
\begin{equation}
\tnccardf{A} = (\tncwa\widehat{\otimes}A_{\chi})^G_{i_a\otimes 1} = \tnccar{A}
\end{equation}
We describe the induced differential and multiplicative structure on $\tnccardf{A}$.

\begin{df}
The twisted nc Cartan differential $d_G^{\chi}$ on $\tnccardf{A}$ is the differential induced by the Kalkman map $\Phi^{\chi}$:
\begin{equation}
\label{defdgnc}
d_G^{\chi} = \Phi^{\chi} ( \nwd\otimes 1 + 1\otimes d) (\Phi^{\chi})^{-1}
\end{equation}
\end{df}
There is a large class of Drinfel'd twists the Cartan differential is insensitive to. A sufficient condition for the equality $d_G^{\chi}=d_G$, as we are going to prove, is that $\chi$ acts as the identity on $\tnccardf{A}$; this is true for example for every $\chi$ depending antisymmetrically by
commuting generators $H_i\in\g$, as it easy to check. For instance the class of Drinfel'd twists relevant for isospectral deformations is of this kind.

\begin{prop}
\label{propnccd}
The differential $d_G^{\chi}$ is the twist of the nc Cartan differential $d_G$ of (\ref{ncdcart}),
$d_G^{\chi}=\chi d_G \chi^{-1}$. In particular, when $\chi$ acts as the identity on $\tnccardf{A}$ we have $d_G^{\chi}=d_G$.
\end{prop}

\proof{ The first statement follows directly from $(\ref{defdgnc})$, using $[\chi,\nwd\otimes 1]=[\chi,1\otimes d]=0$ as operators on $\tncwa\widehat A_{\chi}$; the second part is evident.}

Since so far we discussed Drinfeld twists elements of the type (\ref{twistel}) which satisfies the above conditions, in the following we will use $d_G^{\chi}=d_G$. We can interpret the image of the twisted Kalkman map as a twisted BRST complex, which then restricted to the basic subcomplex gives the twisted nc Cartan model. The twisted nc BRST differential is
\begin{equation}
\label{tncbrstd}
\delta^{\chi}_{BRST} = \Phi^{\chi}(\nwd\otimes 1+1\otimes d)\Phi^{\chi} = \chi(\delta_{BRST})\chi^{-1}
\end{equation}
namely the twist of the nc BRST differential (\ref{ncbrstd}).

The last thing to compute is the multiplicative structure induced in the Cartan complex $(\tnccardf{A},d_G)$; this is determined by $\Phi^{\chi}$ following (\ref{moltcart}). A nice expression is obtained under the following assumption, which is natural if we think of $A$ as the algebra of differential forms.

\begin{prop}
\label{tnccmmolt}
Let us assume $(A,\cdot)$ is graded-commutative and let $(A_{\chi}, \cdot_{\chi})$ be its Drinfeld twist deformation. The multiplication in the Cartan complex $\tnccardf{A}$ is given, for $u_i\otimes\nu_i \in \tnccardf{A}$, by
\begin{equation}
\label{tnccmmolteq}
(u_1\otimes\nu_1)\odot_{\chi}(u_2\otimes\nu_2) = u_1u_2\otimes (-1)^{|\nu_1||\nu_2|}
\cdot_{\chi} \left(\mbox{exp}\{\frac{1}{2} \, B^{ab} i_a\otimes i_b\} (\nu_2\otimes\nu_1) \right)
\end{equation}
\end{prop}

\proof{Since $\tncwa$ and $\ncwa$ have the same algebra structure and we showed that the twisted basic subcomplex is isomorphic to the untwisted one (see Prop \ref{isobastfm}), we can use a formula relating Clifford and wedge products in the odd part of $\tncwa$
\cite{am1}(Lemma $3.1$) $$ \xi_1 \cdot_{Cl} \xi_2 = \wedge \left(\mbox{exp}\{-\frac{1}{2} \, B^{ab} \, i_a\otimes i_b\}(\xi_1\otimes\xi_2)\right) $$
However note that $i_a$ is the untwisted interior derivative, as well as $\wedge$ is the undeformed product. But thanks to Prop \ref{isobastfm} we can nevertheless pass the exponential factor from $\tncwa$ to $A_{\chi}$ on the twisted basic complex as well, so that the remaining part of $\cdot_{\tncwa}$ commutes with $(\Phi^{\chi})_{bas} = (\Phi)_{bas}$. The effect of the braiding on the multiplicative stricture of $(\tncwa\otimes A_{\chi})_{bas}$ is reduced to (\ref{brpr}), so for the moment we have on $u_i\otimes\nu_i \in (\envt{\g}\otimes A_{\chi})^G$ the multiplication rule
$$(u_1\otimes\nu_1)\odot_{\chi}(u_2\otimes\nu_2) = u_1u_2\otimes (-1)^{|\nu_1||\nu_2|}\, \mbox{exp}\{\, \frac{1}{2} \, B^{ab} \, i_a\otimes i_b\}(\nu_2 \cdot_{\chi} \nu_1) $$
In the previous formula the interior product in the exponential are untwisted, since they came from the undeformed Clifford product of the Weil algebra; however using $(\triangle\xi_a)\chi^{-1}=\chi^{-1}(\triangle^{\chi}\xi_a)$ to replace $\cdot_{\chi}$ by the exponential we get the claimed expression in (\ref{tnccmmolteq}) where now the $i_a$ operators are the twisted derivations which act covariantly on $A_{\chi}$.}

Note that for $A_{\chi}=\Omega(\mfdm_{\theta})$ the deformed product $\cdot_{\chi}$ is the nc wedge product $\wedge_{\theta}$ and the induced multiplication on the Cartan model acts like a deformed Clifford product on $\Omega(\mfdm_{\theta})$; moreover note that the arguments $\nu_1$ and $\nu_2$ are switched, as a consequence of the braided product in the Weil model. As in the untwisted case, this ring structure is not compatible with any possible grading in $A$ and gives the twisted nc Cartan model a filtered double complex structure, to be compared with the graded double complex structure of the classical Cartan model. Finally, for $\theta\rightarrow 0$ we get back the product of the untwisted model (\ref{nccm}).

\begin{df}
\label{tnccarmoddef}
The Cartan model for the equivariant cohomology of a twisted nc $\gt$-da $A_{\chi}$ is the cohomology of the complex $(\tnccardf{A},d_G)$:
\begin{equation}
\mathcal{H}^{\chi}_G(A_{\chi}) = \left( \tnccar{A}, \, d_G \right)
\end{equation}
The differential $d_G$ is given in (\ref{ncdcart}); the ring structure $\odot_{\chi}$ of $\tnccardf{A}$ in (\ref{tnccmmolteq}).
\end{df}

\subsection{Examples and reduction to the maximal torus}

We have seen so far that Drinfel'd twists usually generate a sort of 'mild' deformation; many classical results can be adapted to the deformed setting, and non trivial changes appear only when looking at the algebra (for quantities acted) or bialgebra (for the symmetry acting) structures. Therefore we expect that some properties of classical and nc (in the sense of \cite{am1}) equivariant cohomology will still hold in the twisted case, or at least they will have an appropriate corresponding formulation. In this section we show that several classical results can be restated for twisted models.

We begin with the twisted nc basic cohomology ring, or from a geometric point of view the equivariant cohomology of a point; despite its simplicity it plays a crucial role in localization theorems, and by functoriality any equivariant cohomology ring is a module with respect $H_G(\{pt\})$. Of course the Drinfel'd twist does not deform the the algebra $\mathbb{C}$ representing the point, since the $\gt$-da structure is trivial. Let us just apply the definition of the Weil model:
\begin{equation}
\mathcal{H}_G^{\chi}(\mathbb{C}) = H \left( (\tncwa\otimes\mathbb{C})_{bas} , \nwd\otimes 1 \right) = H((\tncwa)_{bas} , \nwd) = (\tncwa)_{bas}
\end{equation}
The last equality is due to $(\nwd)_{|bas}=0$. So the basic cohomology ring for twisted nc equivariant cohomology is $(\tncwa)_{bas}$. The next step is to get a more explicit expression of this ring, and to compare it with the basic rings of nc and classic equivariant cohomology.

For the nc Weil algebra $\ncwa=\env{\gtb}$ the basic subcomplex consists of elements which commute with either even generators ($G$-invariance) and odd generators (horizontality); in other words, it is the center of the super enveloping algebra $\env{\gtb}$. Passing to horizontal generators we are left with $G$-invariant elements of $\env{\g}$, or again the center; this ring is isomorphic, via Duflo
map, to the ring of $G$-invariant polynomials over $\g$. So we have $(\ncwa)_{bas}\simeq (\env{\g})^G \simeq Sym(\g)^G$, and the latter is the basic cohomology ring of classical equivariant cohomology. In $\tncwa$ the actions of $L$ and $i$ are no longer given by commutators with even and odd generators, but by the twisted adjoint action, which is deformed even on single generators; so there is no evident reason why the basic subcomplex should agree with the center. The following shows nevertheless that it is true.

\begin{prop}
\label{basthm}
The basic subcomplex of the twisted nc Weil algebra $\tncwa$ is isomorphic as a ring to $(\ncwa)_{bas}\simeq\env{\g}^G$.
\end{prop}

\proof{ We prove separately the two opposite inclusions; note that the two basic subcomplexes are subalgebras of the same algebra $\ncwa \simeq \tncwa$. Let us start with $X\in (\ncwa)^G_{hor}$; thus $[X,e_a]=[X,\lambda_a]=0$ by (untwisted) $G$-invariance. But
$$L^{\chi}_a(X)=ad_{e_a}^{\chi}(X)= e_aX\lambda_a - \lambda_aXe_a = \lambda_a (e_aX-Xe_a)=0 $$
and similarly
$$i^{\chi}_a(X)=ad_{\xi_a}^{\chi}(X) = \xi_aX\lambda_a - \lambda_aX\xi_a = \lambda_a(\xi_aX-X\xi_a)=0 $$
and so $X\in (\tncwa)_{hor}^G$. On the other hand, take now $Y\in(\tncwa)^G_{hor}$; on Cartan generators the twisted adjoint action still agrees with the commutator, so $[H_i,Y]=0$ and then $[\lambda_a,Y]=0$. But then
$$ad_{e_a}^{\chi}Y=0=e_aY\lambda_a - \lambda_aYe_a=\lambda_a(e_aY-Ye_a)$$
implies the untwisted $ad_{e_a}(Y)=[e_a,Y]=0$; the same for
$$ad_{\xi_a}^{\chi}(Y)=0=\xi_aY\lambda_a - \lambda_aY\xi_a=\lambda_a(\xi_aY-Y\lambda_a)$$
which gives the untwisted $ad_{\xi_a}Y=[\xi_a,Y]=0$. So $Y\in(\ncwa)_{bas}^G$. The linearity follows from the one of operators $L$ and $i$; the ring structures are the same because they descend from the isomorphic algebra structures of $\ncwa \simeq \tncwa$.}

We can then say that classical, nc and twisted nc equivariant cohomologies have the same basic cohomology ring $Sym(\g^{\ast})\simeq \env{\g}^G$ (we identify $\g$ and $\g^{\ast}$ since we are considering quadratic Lie algebras). \\

The next easy example we consider is when the $\g$ action is trivial; algebraically this corresponds to a trivial $\gt$-da structure, i.e. $L$ and $i$ are identically zero. Also in this case the Drinfel'd twist deformation is absent, since its generators act trivially on the algebra. From the Weil model definition we
find
\begin{equation}
\begin{split}
\mathcal{H}_G^{\chi}(A) & = H ( (\tncwa\otimes A)^G_{hor} , \nwd\otimes 1+1\otimes d ) = \\
 & = H ( (\tncwa)^G_{hor}\otimes A, \nwd\otimes 1+1\otimes d ) = (\tncwa)_{hor}^G\otimes H(A) = \\
 & = \env{\g}^G\otimes H(A)
\end{split}
\end{equation}
Thus also in this case the three different models for equivariant cohomology collapse to the same; the only interesting remark is that the $\env{\g}^G$-module structure of $\mathcal{H}^{\chi}_G(A)$ is given by multiplication on the left factor of the tensor product, so that there is no torsion. This is a very special example of a more general class of spaces we are going to mention later for which this phenomenon always takes place; they are called equivariantly formal spaces. \\

We next come to homogeneous spaces. Classically they are defined as the quotient of a (Lie) group $G$ by a left (or right) action of a closed subgroup $K\subset G$; the action is free, so the quotient is a smooth manifold $X=G/K$ on which $G$ still acts transitively, but now with nontrivial isotropy group. We will recall a classical result which leads to a very easy computation of $H_G(G/K)$, and we will extend this idea to twisted nc equivariant cohomology.

There are many interesting homogeneous spaces; we present general results which apply to all of them, but if one prefers to have a specific example in mind, especially in the twisted picture, we suggest the Drinfel'd-twisted sphere $S^4_{\theta}$ acted upon by $\envt{\mathfrak{so}(5)}$ and realized as
the subalgebra of $Fun_{\gamma}(SO(4))$-coinvariants inside $Fun_{\gamma}(SO(5))$ (with $\gamma$ the dual Drinfeld twist of $\chi$, see the discussion after Thm (\ref{covdef}) and \cite{maj}).

In the classical setting, we consider commuting actions of two Lie groups $K_1$ and $K_2$. If we define $G=K_1\times K_2$ its Weil algebra decomposes in $W_{\g} = W_{\mathfrak{k}_1}\otimes W_{\mathfrak{k}_2}$ with $[\mathfrak{k}_1,\mathfrak{k}_2]=0$ by commutativity of the actions. Then every $\gt$-da algebra $A$ can be thought separately as a $\tilde{\kk}_{1,2}$-da and the basic subcomplex can be factorized in both ways
\begin{equation}
A_{bas \, \g} = (A_{bas \, \kk_1})_{bas \, \kk_2} = (A_{bas \, \kk_2})_{bas \, \kk_1}
\end{equation}

\begin{prop}
\label{commact}
Under the previous assumptions and notations, if $A$ is also both $\kk_1$ and $\kk_2$ locally free we have
\begin{equation}
H_G(A) = H_{K_1}(A_{bas \, \kk_2}) = H_{K_2}(A_{bas \, \kk_1})
\end{equation}
\end{prop}

We simply apply the definition of the Weil model and make use of the commutativity between the two locally free $K_1$ and $K_2$ actions; see e.g. \cite{gs} for the complete proof. This easy fact is very useful for computing equivariant cohomology of homogeneous spaces $H_G(G/K)$. Indeed take on $G$ the two free actions of $K$ and $G$ itself by multiplication; we make them commute by considering $K$ acting from the right and $G$ from the left, or vice versa. The hypothesis of Prop \ref{commact} are satisfied, so we quickly have
\begin{equation}
\label{coomhom}
H_G(G/K) = H_K( G \backslash G) = H_K (\{pt\}) = Sym(\kk^{\ast})^K
\end{equation}
We want to find a similar result for twisted $\gt$-da's. The definition of commuting actions makes perfectly sense in the twisted setting: we require that the two twisted $\tilde{\kk}_{1,2}$-da structures commute. This is an easy consequence of the commutation of the actions on classical algebras, provided
the generators of the twists commute with each other (for example using a unique abelian twist for both algebras, which is the most common situation). The assumption of the local freeness of the action is a bit trickier; we need a good definition of this notion for twisted nc algebras. We refer to the next subsection
for this point; we use that if $A$ is a locally free $\gt$-da then $A_{\chi}$ is a locally free twisted $\gt$-da. So we can apply Prop \ref{commact} also to Drinfeld twist deformations of homogeneous spaces, since all the hypotheses are still satisfied. The appropriate statement involves Drinfeld twists on function algebras over classical groups; this is a dual Drinfel'd twist (see discussion after Thm \ref{covdef}) which deforms the algebra rather than the coalgebra structure of an Hopf algebra. 
We denote by $\gamma:Fun(G)\otimes Fun(G)\rightarrow\mathbb{C}$  the generator of the dual Drinfel'd twist on $Fun(G)$ (dual with respect the $\chi$ twist on $\env{\g}$), which satisfies $\langle \chi,\gamma\rangle =1$ where the brackets come from the duality between $\env{\g}$ and $Fun(G)$. Then the restatement of (\ref{coomhom}) is
\begin{equation}
\label{nccoomhom}
\mathcal{H}^{\chi}_G( (Fun_{\gamma}(G))^{co K} ) = \mathcal{H}^{\chi}_K( (Fun_{\gamma}(G))^{co G} ) = \mathcal{H}^{\chi}_K (\mathbb{C}) =
\env{\kk}^K
\end{equation}
As an explicit example, we can apply (\ref{nccoomhom}) to nc spheres $S^n_{\theta}$. For simplicity let us consider $S^4_{\theta}$; it can be constructed as a toric isospectral deformation of the classical sphere $S^4$ twisting the $\mathbb{T}^2$ symmetry acting on it. Equivalently, to stress the fact that it is a homogeneous space, we can think of it as the $Fun_{\gamma}(SO(4))$-coinvariant subalgebra of
$Fun_{\gamma}(SO(5))$. On $S^4_{\theta}$ we have the action of the twisted symmetry $\envt{so(5)}$; the action of course is not free since the twisted Hopf subalgebra $\envt{so(4)}$ acts trivially. The equivariant cohomology of this twisted action is defined using the twisted Weil (or Cartan)
models introduced in the previous section, and it may be computed using (\ref{nccoomhom}). We find
\begin{equation*}
\mathcal{H}_{\mathfrak{so}(5)}^{\chi}(S^4_{\theta}) = \envt{\mathfrak{so}(4)}^{SO(4)} = \env{\mathfrak{so}(4)}^{SO(4)}
\simeq Sym(\mathfrak{so}(4))^{SO(4)} \simeq Sym(\ttt^2)^W
\end{equation*}
where the last equality is given by Chevalley's theorem $Sym(\g)^G\simeq Sym(\ttt)^W$ for $W$ the Weyl group. \\

We now study the reduction of twisted nc equivariant cohomology to the maximal torus $T\subset G$.
The two main ingredients in the algebraic proof of the isomorphism $H_G(X)=H_T(X)^W$ ($W$ denotes the Weyl group of $T$) for classical equivariant cohomology are the functoriality of $H_G(X)$ with respect to group reduction $P\subset G$, and spectral sequence arguments.

In order to reproduce a similar result and proof for the nc (and then twisted) case we first need to work out the functorial properties of $\mathcal{H}_G(A)$; since in both nc and twisted cases Weil and Cartan models are built using the Lie algebra $\g$, contrary to the classical case which makes use of the dual $\g^{\ast}$, it is not obvious that for every subgroup $P\subset G$ we have a morphism of Cartan complexes $\subnccar{P}{A}\rightarrow\nccar{A}$. The existence of such a morphism is guaranteed for the specific choice $P=N(T)$, the normalizer of the maximal torus, by a generalized Harish-Chandra projection map \cite{am2}. Once we have such morphism, the rest of the proof follows quite easily.

We start with a quick review of the classical reduction, referring to \cite{gs} for the full details. The Cartan complex $C_G(A)$ may be seen as a double $\mathbb{Z}$-graded complex $C^{p,q}(A) = (Sym^p(\g^{\ast})\otimes A^{q-p})^G$ with differentials $\delta_1=-v^a\otimes i_a$ and $\delta_2=1\otimes d$ of grading $(1,0)$ and $(0,1)$ respectively; the cohomology of the total complex with respect to $d_G=\delta_1 + \delta_2$ is the classical equivariant cohomology. This gives the usual setting to construct a spectral sequence converging to $H_G(A)$  with $E_1^{p,q}$ term (for $G$ compact and connected) given by $Sym^p(\g^{\ast})\otimes H^{q-p}(A)$. We can get the desired isomorphism $H_G(X)\cong H_T(X)^W$ by looking at a different spectral sequence having the same $E_1$ term. For each closed subgroup $P\subset G$ we get a morphism between Cartan complexes $C_G(A)\rightarrow C_P(A)$ and hence between $E_1$ terms; whenever $P$ is such that $Sym(\g^{\ast})^G\cong Sym(\p^{\ast})^P$ we have an isomorphism at the $E_1$ step compatible with the differentials, thus it descend to every following step and in particular $H_G(A)\cong H_P(A)$. We can use this result with $P=N(T)$, the normalizer of the maximal torus.

\begin{thm}
\label{redmaxt}
Let $G$ be a compact connected Lie group and $A$ a $\gt$-da. There is a ring homomorphism $H_G(A)\cong H_T(A)^W$ where $T\subset G$ is the maximal torus in $G$ and $W$ its Weil group $N(T)/T$.
\end{thm}
\proof{The Weil group $W=P/T=N(T)/T$ is finite, thus $\p\cong\ttt$ and $Sym(\p^{\ast})^P\cong Sym(\ttt^{\ast})^P\cong Sym(\ttt^{\ast})^W$ since $T$ acts trivially on itself. Then by Chevalley's theorem $Sym(\g^{\ast})^G\cong Sym(\ttt^{\ast})^W$, so as discussed before $H_G(A)\cong H_{N(T)}(A)$. To conclude we have to prove that $H_{N(T)}(A)\cong H_T(A)^W$; the inclusion $T\hookrightarrow P=N(T)$ induces a morphism $Sym(\p^{\ast})\otimes A\rightarrow Sym(\ttt^{\ast})\otimes A$ and by taking the $P$-invariant subcomplexes we get a morphism $C_P(A)\rightarrow C_T(A)^W$ and so on at each stage of the spectral sequences. In particular we obtain a morphism between equivariant cohomologies $H_P(A)\rightarrow H_T(A)^W$; but note that at the $E_1$ step the morphism is indeed an isomorphism, since $Sym(\p^{\ast})^P\cong Sym(\ttt^{\ast})^W$, so the previous morphism between cohomologies is an isomorphism as well.}

This result allows us to reduce the computation of classical equivariant cohomology for generic compact Lie groups $G$ to abelian groups. Another important feature of $H_G(X)$ is its $Sym(\g^{\ast})^G$-module structure, with the torsion part playing a central role in localization theorems. We proved that the $E_1$ term of the spectral sequence converging to $H_G(X)$ is $Sym(\g^{\ast})^G\otimes H(A)$; at this stage the module structure is simply given by left multiplication, so $E_1$ is a free $Sym(\g^{\ast})^G$-module. This already implies that if $H(A)$ is finite dimensional, the equivariant cohomology ring $H_G(A)$ is finitely generated as $Sym(\g^{\ast})^G$-module. When the spectral sequence collapses at this stage, the algebra $A$ is called equivariantly formal. The definition comes from \cite{gkp} (using the language of $G$-spaces $X$ rather then $\gt$-da's $A$), where sufficient conditions for the collapsing are studied. In this case since $E_{\infty}\cong E_1$ we have that $H_G(A)$ is a free $Sym(\g^{\ast})^G$-module. We can also express the ordinary cohomology in terms of equivariant cohomology by tensoring the $E_1$ term by the trivial $Sym(\g)^{\ast}$-module $\mathbb{C}$, obtaining $H(A)=\mathbb{C}\otimes_{Sym(\g^{\ast})}H_G(A)$.

We now come to nc equivariant cohomology. Given a closed subgroup $P\subset G$ we have a Lie algebra homomorphism $\p\rightarrow\g$ which may be lifted to the enveloping algebras and nc Weil algebras, but in general does not intertwine the differentials and most unpleasantly goes in the opposite direction to the one in which we are interested in order to reduce equivariant cohomology. We have to look for a $\pt$-da (or at least $\pt$-ds, i.e. $\pt$-differential space) homomorphism $\ncwa\rightarrow\mathcal{W}_{\p}$ which then may be used to get a morphism between the nc Cartan complexes $\nccar{A}\rightarrow\subnccar{P}{A}$.
This homomorphism can be constructed for a very special choice of the subgroup $P$, namely for $P=N(T)$, which is exactly the case we need. We refer to \cite{am2}(Section $7$) for the details of the construction. It is shown that for a quadratic Lie algebra $\g$ with quadratic subalgebra $\p$ and orthogonal
complement $\p^{\bot}$ it is possible to define a 'generalized' Harish-Chandra projection
$k_{\mathcal{W}}:\ncwa\rightarrow\mathcal{W}_{\p}$ which is a $\pt$-ds homomorphism and becomes a $\pt$-da homomorphism between the basic subcomplexes $\env{\g}^G\rightarrow\env{\p}^P$. Moreover this construction reduces to the classical Harish-Chandra map up to $\p$-chain homotopy \cite{am2}(Thm $7.2$)
and then looking at the basic subcomplexes (where the differential is zero) we find the commutative diagram of $\pt$-da's \cite{am2}(Thm $7.3$)
\begin{equation}
\label{diagrhc}
\xymatrix{
Sym(\g)^G\ar[r]\ar[d]_{k_{Sym}}& \env{\g}^G\ar[d]^{(k_{\mathcal{W}})_{|bas}} \\
Sym(\p)^P\ar[r]               & \env{\p}^P}
\end{equation}
where horizontal maps are Duflo algebra isomorphism. For $P=N(T)$ by Chevalley's theorem the map $k_{Sym}:Sym(\g)^G\rightarrow Sym(\ttt)^W$ is an algebra isomorphism as well. This is the morphism we need
to prove the reduction of nc equivariant cohomology. We note that this reduction Thm, even if not explicitly stated, is already contained in \cite{am1} when the authors prove the ring isomorphism $H_G(A)\cong\mathcal{H}_G(A)$ induced by the quantization map $Q_{\g}:W_{\g}\rightarrow\ncwa$. We prefer to give here a direct proof based on morphisms between Cartan complexes and spectral sequences since this approach will be generalized to our twisted nc equivariant cohomology.
\begin{thm}
\label{ncredmaxt}
The ring isomorphism of Thm(\ref{redmaxt}) holds also between nc equivariant cohomology rings; for every nc $\gt$-da $A$ and compact connected Lie group $G$ the reduction reads $\mathcal{H}_G(A)\cong \mathcal{H}_T(A)^W$.
\end{thm}
\proof{As for the classical reduction, the proof is based on the presence of a morphism between Cartan complexes and a comparison between the two associated spectral sequences. The setting is now the following: the nc Cartan model $\nccar{A}=(\env{\g}\otimes A)^G$ is looked at as a double filtered differential complex. On one side we have the standard increasing filtration of the enveloping algebra $\env{\g}_{(0)}\subset\env{\g}_{(1)}\subset\env{\g}_{(2)}\ldots$; on the other side, supposing $A$ is a finitely generated graded algebra, we have an increasing filtration $A_{(p)}=\oplus_{i\leq p}A^i$; note that this double filtration on $\nccar{A}$ is compatible with the ring structure (\ref{nccm}) (contrary to the grading of $A$, which is not compatible with the induced product on $\nccar{A}$). The operators
$$\delta_1=\Phi(d^{\ncwa}\otimes 1)\Phi^{-1}=- \, \frac{1}{2} (u^a_{(L)}+u^a_{(R)})\otimes i_a + \frac{1}{24}f^{abc}\otimes i_ai_bi_c$$ and $$\delta_2=\Phi(1\otimes d)\Phi^{-1}=1\otimes d$$
square to zero (since their counterpart on the Weil complex do), and then anti-commute since their sum is the nc Cartan differential $d_G$; they are the differentials of the double complex, with filtration degree respectively $(1,0)$ and $(0,1)$. The cohomology of the total complex with respect to
$d_G=\delta_1+\delta_2$ is the nc equivariant cohomology ring $\mathcal{H}_G(A)$; the filtration of $\nccar{A}$ induces a filtration on the cohomology. We can compute its graded associated module $Gr(\mathcal{H}_G(A))$ by a spectral sequence with $E_0$ term given by the graded associated module of the nc Cartan model $Gr(\nccar{A})=C_G(A)$; this is the spectral sequence we already introduced before. Note that the differentials $\delta_1$ and $\delta_2$ map to the ordinary differentials of the Cartan complex $-\frac{1}{2}v^a\otimes i_a$ and $1\otimes d$. Now let us consider the inclusion $P=N(T)\subset G$ and the Harish-Chandra projection map $k_{\mathcal{W}}:\ncwa\rightarrow\mathcal{W}_{\p}$. This induces a $\pt$-ds morphism between the Weil complexes $(\ncwa\widehat{\otimes}A)_{bas}\rightarrow (\mathcal{W}_{\p}\widehat{\otimes}A)_{bas}$ and by Kalkman map a $\pt$-ds morphism between nc Cartan models $\nccar{A}\rightarrow\subnccar{P}{A}$ compatible with the filtrations; commuting with differentials, it also lifts to cohomology giving a morphism of filtered rings $\mathcal{H}_G(A)\rightarrow\mathcal{H}_P(A)$. By going to the graded associated modules and computing the $E_1$ term of the spectral sequence we get a $\pt$-ds morphism $Sym(\g)^G\otimes H(A) \rightarrow Sym(\ttt)^W\otimes H(A)$ (see (\ref{diagrhc}) and \cite{am2}(Thm $7.3$)). Now this is a $\pt$-da isomorphism, and it induces $\pt$-da isomorphisms at every further step of the spectral sequence. The isomorphism between $Gr(\mathcal{H}_G(A))$ and $Gr(\mathcal{H}_P(A)$ implies that the morphism $\mathcal{H}_G(A)\rightarrow\mathcal{H}_P(A)$ introduced before is in fact a ring isomorphism. As in the classical case, the last step is to show $\mathcal{H}_P(A)\cong\mathcal{H}_T(A)^W$; this easily follows from the morphism $\subnccar{P}{A}\rightarrow\subnccar{T}{A}$ (note that $\p\cong\ttt$ so the previous morphism is just group action reduction) and a completely similar spectral sequence argument.}

We finally note that another equivalent proof of Thm \ref{ncredmaxt} may be obtained by a different construction of the morphism $\nccar{A}\rightarrow\subnccar{P}{A}$ via a diagram
\begin{equation}
(\env{\p}\otimes A)^P\longrightarrow ((\env{\g}\otimes Cl(\p^{\bot}))\otimes A)^P \longleftarrow (\env{\g}\otimes A)^G
\end{equation}
Considering the spectral sequence associated to these three Cartan models (the cohomology of the middle complex is a sort of 'relative' equivariant cohomology $\mathcal{H}_{G,P}(A)$ of $G$ with respect to $P$, see \cite{am2}(Section $6$)) it is possible to prove an isomorphism between the image of the left and right $E_1$ terms inside the $E_1$ term of the middle complex \cite{am2}(Thm $6.4$). This isomorphism is referred as a version of Vogan's conjecture for quadratic Lie algebras. 

We finally consider twisted nc equivariant cohomology. It is a natural question to ask if our model satisfies a reduction property as well; an easy nevertheless crucial fact is that Drinfeld twists act trivially on abelian symmetries. This will allow us to basically use the same proof of Thm \ref{ncredmaxt};
moreover for the same reason when restricted to the maximal torus $T$, twisted nc equivariant cohomology $\mathcal{H}^{\chi}_T(A_{\chi})$ agrees with $\mathcal{H}_T(A_{\chi})$.
\begin{thm}
\label{tncredmaxt}
Let $G$ be a compact connected Lie group, and $A_{\chi}$ a twisted $\gt$-da. There is a ring homomorphism
$\mathcal{H}^{\chi}_G(A_{\chi})\cong \mathcal{H}^{\chi}_T(A_{\chi})^W$ where $T\subset G$ is the maximal torus in $G$ and $W$ its Weil group $N(T)/T$.
\end{thm}
\proof{We can use the generalized Harish-Chandra projection also for twisted nc Weil algebras, since for $P=N(T)$ as $\pt$-da's $\ncwa\cong\tncwa$. The twisted nc Cartan model $\tnccardf{A}$ is a double filtered differential complex similarly to $\nccar{A}$, and we can consider the spectral sequence constructed from its graded associated module. At the $E_1$ step as usual we are left with the basic part of $Gr(\tncwa)$ tensored with $H(A_{\chi})$; since $(\tncwa)_{|bas}\cong (\ncwa)_{|bas}$ (see Thm \ref{basthm}) any effect of the twist is now present only in the cohomology of $A_{\chi}$. Then the isomorphism between the $E_1$ terms of $\tnccardf{A}$ and $\subtnccardf{P}{A}$ follows as in the proof of Thm \ref{ncredmaxt}. The same happens for the last part of the proof, when going from $P=N(T)$ to $T$.}

This result shows one more time that deformations coming from Drinfel'd twists do not affect much of the classical setting. The definition of a twisted nc equivariant cohomology is needed when dealing with algebras which carry a twisted action of a symmetry, and this is exactly what happens for covariant actions of Drinfel'd twisted Hopf algebras. However the possibility to reduce the cohomology to the maximal torus part leaves the only contribution coming from the Drinfeld twist in the deformed ring structure of $\mathcal{H}^{\chi}(A_{\chi})$, while the vector space and $Sym(\g)^G$-module structures are undeformed.

The positive part of this quite classical behaviour is that for what concerns this class of deformations, a lot of techniques of equivariant cohomology may be lifted with an appropriate and careful rephrasing to the nc setting. On the contrary, if we are interested in purely new phenomena which do not admit a classical
counterpart, it seems we have to enlarge the class of deformations considered, either by taking Drinfel'd twists $\chi$ which do not satisfy the $2$-cocycle condition or moving to other classes of deformations. To this end we present in the next subsection a sketch of a general strategy to define Weil models for equivariant cohomology of more general class of deformations.

\subsection{Models for generic deformations}

We briefly outline in this last subsection a general approach towards a definition of algebraic models for the equivariant cohomology of deformed $\gt$-da's. This is the relevant formalism for nc spaces which carry a covariant action of some deformed symmetry.

Indeed we can reinterpret the above described models for twisted nc equivariant cohomology as a particular example of a more general construction. We present this general construction by focusing on five steps. We have two ideas in mind: first, we can apply this plan to different classes of deformations, for example Drinfel'd-Jimbo quantum enveloping algebras and their covariant actions, and study the associated nc equivariant cohomology \cite{lccp}. On the other hand we feel that this general approach may cast some light on the twisted models themselves, in particular on the role played by our twisted nc Weil algebra and its universality. For example it turns out that a simpler Weil algebra can be used to define the cohomology, leading to a possible easier expression of the models. The full details on this new formulation of the twisted models, as well as the proof of the results we claim here will appear in \cite{lu2}. 

We summarize the strategy by listing five sequential steps; we then discuss more carefully each of them, and we make some further comment about how they fit with our definition of twisted nc equivariant cohomology.

\begin{enumerate} 
\item Choose the relevant category of Hopf-module algebras. This amounts to choose the deformed $\gt$-da structure, i.e. the deformation of the symmetry and/or of the nc space acted.
\item Give a suitable definition of locally free action in the category. Equivalently, characterize algebraic connections on the algebras of the category considered.
\item Find the universal locally free algebra of the category; this object $\mathcal{W}^{\prime}$ will be interpreted as the deformed Weil algebra associated to the choosen class of deformations.
\item For each algebra $A$ in the category define the Weil model for equivariant cohomology as the cohomology of the basic subcomplex $(\mathcal{W}^{\prime}\otimes A)_{bas}$.
\item For the Cartan model, consider a deformation of the Kalkman map compatible with the deformation of the category of Hopf-module algebras.
\end{enumerate}

The first point summarizes the fact that to consider covariant actions of symmetries is equivalent to work in the category of Hopf-module algebras. Starting with a deformed nc algebra $A_{\theta}$ we realize a covariant action of a classical symmetry $\g$ by realizing a deformed $\gt$-da structure on $A_{\theta}$, i.e. by fixing the relevant category of Hopf module algebras to which $A_{\theta}$ belongs. 

Once we fix the category, so we have a compatible deformation of symmetries and spaces, we need to distinguish locally free actions. We know how equivariant cohomology is defined for locally free actions, and we want to reduce every other case to a locally free action. Classically a $\gt$-da $A$ carries a locally free action if it admits an algebraic connection; for $\g$ quadratic this is equivalent to a $\gt$-da morphism $\vartheta:Sym(\gt)\rightarrow A^1$ (if $A$ is graded we want the image to have degree one) \cite{am2}. Thus an algebraic connetion is a morphism in the category between the symmetric $\gt$-da and the algebra considered. This can be generalized to arbitrarily deformed $\gt$-da's; we only need to consider the deformed symmetric $\gt$-da and ask for the connection to be a morphism in the deformed category.  

As next step, in analogy with the classical definition, we interpret the universal locally free object in the category of deformed $\gt$-da's as a deformed Weil algebra. Looking at the definition of algebraic connection a natural candidate is the deformed symmetric $\gt$-da itself, endowed with a Koszul differential that ensures aciclicity. In some sense, in order to encode the deformation of the category, our definition of algebraic connection is already given at the level of the induced Chern-Weil morphism, so that it comes directly associated with a Weil algebra. Note that as in the classical case, any algebra in the category which is $\gt$-homotopic with $\mathcal{W}^{\prime}$ (we call such algebras of Weil-type using the terminology of \cite{am2}, or  $\mathcal{W}^{\ast}$-modules following \cite{gs}) can be used in place of $\mathcal{W}^{\prime}$ to define equivariant cohomology. 

A Weil model for equivariant cohomology is then defined by considering the tensor product in the category of deformed $\gt$-da's between the deformed Weil algebra $\mathcal{W}^{\prime}$ and the algebra we want to take cohomology. Note that this tensor product is in general braided, depending on the quasitriangular structure of the deformation of $\env{\gt}$. The notion of basic subcomplex still makes sense, since the deformed $\gt$-da structure provides deformed Lie and interior derivatives acting on the algebras of the category.

Finally, if one wants to pass from the deformed Weil model to a deformed Cartan model, a suitable Kalkman map has to be constructed; following \cite{kal}, we interpret the image of this Kalkman map as a deformed BRST model, while its restriction to the basic subcomplex defines the deformed Cartan model. \\

We quickly show how this strategy reflects what we actually have done dealing with Drinfel'd twist deformations. To this class of deformations corresponds the category of $\envt{\gt}$-module algebras; as shown by Drinfel'd \cite{dr1}\cite{dr2} this category is equivalent to the undeformed one, and we have the explicit tensor functor which realizes the equivalence. Following the claimed plan, we could take as twisted nc Weil algebra the twisted symmetric $\gt$-da, which can be defined as the quotient of the tensor algebra of $\gt$ by the braided-symmetric relations $a\otimes b - \Psi(a,b)=0$ ($\Psi$ is the braiding morphism of the category, induced by $\chi$). What we have done in the present paper is a bit different; we started with the Weil algebra of \cite{am1} and deformed it by a Drinfled twist. The reason is that we realized this general strategy only recently. Our claim is that the same twisted models may be defined in an equivalent (and maybe simpler, expecially at the level of Cartan complex) way by using the twisted symmetric algebra as deformed Weil algebra, and that our $\tncwa$ is actually of Weil-type and twisted $\gt$-homotopic to the 'real' Weil algebra. We plan to discuss these topics in \cite{lu2} and to apply this five-steps construction to Drinfel'd-Jimbo deformations in \cite{lccp}.



\addcontentsline{toc}{section}{References}
\frenchspacing
\begin{thebibliography}{9}

\bibitem[AM00]{am1} A. Aleskeev, E. Meinrenken: \emph{The non-commutative
Weil algebra}, Invent. Math. 139 (2000), 135-172

\bibitem[AM05]{am2} \bysame : \emph{Lie theory and the Chern-Weil homomorphism},
Ann. Scient. Ec. Norm. Sup. 38 n. 4 (2005), 303-338

\bibitem[Car50]{car} H. Cartan: \emph{Notions d'alg\`{e}bre diff\'{e}rentielle; application aux groupes de Lie at aux vari\'{e}t\'{e}s o\`{u} op\`{e}re un groupe de Lie}, Colloque de Topologie, {C.B.R.M.}, Bruxelles (1950), 15-27 

\bibitem[CVD02]{cdv} A. Connes, M. Dubois-Violette: \emph{Noncommutative finite-dimensional manifolds $I$. Spherical manifolds and related examples}, Comm. Math. Phys. 230 n. 3 (2002), 539-579

\bibitem[Cir]{lu2} L. Cirio. In preparation.

\bibitem[CL01]{cl} A. Connes, G. Landi: \emph{ Noncommutative manifolds, the instanton algebra
and isospectral deformations}, Comm. Math. Phys. 221 n. 1 (2001), 141-159

\bibitem[CP]{lccp} L. Cirio, C. Pagani. In preparation.

\bibitem[DG97]{qla} G.W. Delius, M.D. Gould: \emph{Quantum Lie algebras, their existence, uniqueness and $q$-antisymmetry}, Comm. Math. Phys. 185 n. 3 (1997), 709-722

\bibitem[Dri90a]{dr1} V.G. Drinfel'd: \emph{On almost cocommutative Hopf algebras}, Leningrad Math. J. 1 (1990), 321-342

\bibitem[Dri90b]{dr2} \bysame : \emph{Quasi-Hopf algebras}, Leningrad Math. J. 1 (1990), 1419-1457

\bibitem[GKR98]{gkp} M. Goresky, R. Kottwitz, R.MacPherson: \emph{Equivariant cohomology, Koszul duality, and the localization theorem}, Invent. Math. 131 (1998), 25-83

\bibitem[GS99]{gs} V. Guillemin, S. Sternberg: \emph{Supersimmetry and equivariant de Rham theory}, Springer-Verlag, (1999)

\bibitem[Kal93]{kal} J. Kalkman: \emph{BRST model for equivariant cohomology and representatives for the equivariant Thom class}, Comm. Math. Phys. 153 n. 3 (1993), 447-463

\bibitem[Kas95]{kas} C. Kassel: \emph{Quantum groups}, Graduate Texts in Mathematics vol 115, Springer-Verlag, (1995)

\bibitem[Maj94]{maj} S. Majid: \emph{Foundations of quantum group theory}, Cambridge University Press, (1994)

\bibitem[Res90]{res} N. Reshetikhin: \emph{Multiparameter quantum groups and twisted quasitriangular Hopf algebras}, Lett. Math. Phys. 20 n. 4 (1990), 331-335

\bibitem[Rie93]{rie} M. Rieffel: \emph{Deformation quantization for actions of $\mathbb{R}^d$}, Mem. Amer. Math. Soc. vol 506, (1993)

\bibitem[SS93]{ss} S. Shnider, S. Sternberg: \emph{Quantum groups - from coalgebras to Drinfel'd algebras}, Graduate texts in Mathematical Physics vol. 2, Cambridge International Press, (1993)

\end{thebibliography}
\end{document}